\documentstyle[amsfonts,amssymb,amsfonts,amssymb,graphicx,mathrsfs,color,ulem]{article}

\newcommand{\R}{{\mathbb R}}
\newcommand{\C}{{\mathbb C}}

\newcommand{\dsp}{\displaystyle}

\newcommand{\ol}{\overline}
\newcommand{\ndt}{\noindent}
\newcommand{\qed}{\hfill $\square$}

\newtheorem{theorem}{Theorem}[section]
\newtheorem{lemma}[theorem]{Lemma}

\newtheorem{corollary}[theorem]{Corollary}

\title{\bf 
The  self-interacting Dirac fields in FLRW spacetime
}

\author{{\bf Anahit Galstian and Karen Yagdjian} }

\begin{document}

\date{}
\maketitle
\thispagestyle{empty}
\vspace{-0.3cm}

\begin{center}
{\it School of Mathematical and Statistical Sciences,\\
University of Texas RGV,\\
1201 W.~University Drive,  
Edinburg, TX 78539,
USA }
\end{center}
\medskip

\addtocounter{section}{-1}
\renewcommand{\theequation}{\thesection.\arabic{equation}}
\setcounter{equation}{0}
\pagenumbering{arabic}
\setcounter{page}{1}
\thispagestyle{empty}

\hspace{2cm}\begin{abstract}
\begin{small}
We prove the existence of global in time  solution with the small initial data for the semilinear equation of the spin-$\frac{1}{2}$ particles in the Friedmann-Lema\^itre-Robertson-Walker spacetime. Moreover, we also prove that if the initial function satisfies the Lochak-Majorana condition, then the global solution exists for arbitrary large initial value. The solution scatters to free solution for large time. The mass term is assumed to be complex-valued. The conditions on the imaginary part of mass are discussed by proving nonexistence of the global solutions if certain relation between scale function and the mass are fulfilled.  
\medskip

\end{small}
\end{abstract}

\renewcommand{\theequation}{\thesection.\arabic{equation}}

\section{Introduction}
\setcounter{equation}{0}

The behavior of particles obeying the covariant Dirac equation in the Friedmann-Lema\^itre-Robertson-Walker (FLRW)  spacetime has given rise to many investigations { (\cite{Barut-D,Fock,Huang_Justin,Oliveira,Parker,Schrodinger,Zecca})}. On the other hand, the solution of the semilinear Dirac equation describing self-interacting Dirac fields in general FLRW spacetime according to our best knowledge, haven't been discussed in literature.
In this article we study   solutions of the semilinear Dirac  equation in the curved spacetime of the FLRW models of cosmology. In particular, we focus on  the relationship between the mass term, scale factor,  nonlinear term, and initial function, which provides   a global in time existence  or an estimate on the lifespan  of the solution  of the Dirac equation in the expanding universe.

The review of some results { (\cite{Federico-Piero,Candy,Piero,Shuji})} on  small and large amplitude global solutions of the semilinear Dirac equation in Minkowski
space  is given in  \cite{ArX2021b}.

The metric tensor of the spatially flat FLRW spacetime in Cartesian coordinates  is 
\begin{equation}
\label{metric}
(g_{\mu \nu })=    \left (
   \begin{array}{ccccc}
 1& 0& 0   & 0 \\
   0& -a^2(t) &  0 & 0 \\ 
 0 & 0 &  -a^2(t)   & 0 \\
 0& 0& 0   &  -a^2(t) \\
   \end{array} \right),\quad \mu ,\nu =0,1,2,3,
\end{equation}
where the scale factor    $a(t)=a_0t^{\ell}  $, $\ell \in {\mathbb R}$,  $ t>0$,  and $x \in {\mathbb R}^3$, $x_0=t $. If $\ell<0$ the spacetime is contracting. In the case of $\ell>1$ the expansion is accelerating (with horizon), while for $0<\ell<1$ the expansion is decelerating. 
In particular, the Einstein-de~Sitter space with the scale factor   $a(t)=t^{2/3}  $  is modeling the expanding  matter dominated universe, if $a(t)=t^{1/2}  $ is radiation dominated universe (see, e.g., \cite{Moller,Ohanian-Ruffini}). The scale factor   $a(t)=t   $ describes the Milne model \cite{Sean Carroll,Gron-Hervik, Schrodinger}. The curvature of the space  with   $a(t)=a_0t^{\ell}  $ is  $R(t)=- 6 \ell \left(1-2 \ell  \right) t^{-2}   $. The Dirac equation in the   space with metric (\ref{metric}) is 
\begin{equation}
\label{DE}
   \dsp 
\left(  i {\gamma }^0    \partial_0   +i \frac{1}{a(t)} \sum_{k=1,2,3}{\gamma }^k  \partial_{ k} +i \frac{3\dot a(t)}{2a(t)}    {\gamma }^0     -m{\mathbb I}_4 \right)\Psi=f \,,
\end{equation}
(see, e.g., \cite{Barut-D}), where { $\partial_0=\partial /\partial t $, $\partial_k=\partial /\partial x_k $ and} $\frac{ \dot a(t)}{ a(t)}$ is the Hubble parameter { (see, e.g.,  \cite[Ch.8]{Sean Carroll},\cite[Sec.11.4]{Gron-Hervik})}, while 
\[
 \gamma ^0= \left (
   \begin{array}{ccccc}
   {\mathbb I}_2& 0   \\
   0& -{\mathbb I}_2   \\ 
   \end{array}
   \right),\quad 
\gamma ^k= \left (
   \begin{array}{ccccc}
  {\mathbb O}_2& \sigma ^k   \\
  -\sigma ^k &  {\mathbb O}_2  \\  
   \end{array}
   \right),\quad k=1,2,3,\quad
 \gamma ^5:=-i\gamma ^0\gamma ^1\gamma ^2\gamma ^3=\left(
\begin{array}{cccc}
 {\mathbb O}_2 & -{\mathbb I}_2 \\
 -{\mathbb I}_2 & {\mathbb O}_2 \\
\end{array}
\right) \nonumber .
 \]
Here $\sigma ^k $ are Pauli matrices 
\begin{eqnarray*}
&  &
\sigma ^1= \left (
   \begin{array}{ccccc}
  0& 1   \\
  1& 0  \\  
   \end{array}
   \right), \quad
\sigma ^2= \left (
   \begin{array}{ccccc}
  0& -i   \\
  i& 0  \\  
   \end{array}
   \right),\quad
\sigma ^3= \left (
   \begin{array}{ccccc}
  1& 0   \\
  0&-1 \\  
   \end{array}
   \right)\,,
\end{eqnarray*}
and  ${\mathbb I}_n $, ${\mathbb O}_n $ denote the $n\times n$ identity and zero matrices, respectively.

Consider the Dirac equation
\begin{eqnarray*}
   \dsp 
\left(  i {\gamma }^0    \partial_0   +i \frac{1}{a(t)} \sum_{k=1,2,3}{\gamma }^k  \partial_{k}  +i \frac{3\dot a(t)}{2a(t)}    {\gamma }^0     -m{\mathbb I}_4 +{\gamma }^0V(x,t) \right)\Psi=0  \,,
\end{eqnarray*}
where $V$ is the matrix-valued potential $  V(x,t): {\mathbb R}^3\times (0,\infty) \longrightarrow M_4({\mathbb C}) $. 
  The last equation can be written in the equivalent form as following,
\begin{eqnarray*}
   \dsp 
\left(   \partial_0 +\frac{1}{a(t)}\sum_{j=1,2,3}\alpha ^j  \partial_{j} + \frac{3\dot a(t)}{2a(t)}   {\mathbb I}_4    +im{\gamma }^0  
-iV (x,t)\right)\Psi=0  \,,
\end{eqnarray*}
where ${\alpha }^j={\gamma }^0{\gamma }^j $ are self-adjoint matrices,  $({\alpha }^j)^*={\alpha }^j$. 
The semilinear Dirac equation, which describes the self-interacting field, can be obtained by introducing into the  equation the   nonlinear term. { We characterize  that term by the next  condition. }

\noindent
{\bf Condition (${\mathcal L}$)} {\it The function $F=F(\psi ) \in C^3({\mathbb C}^4;{\mathbb C}^4)$   is Lipschitz continuous with exponent $\alpha  $ in the space $H_{(s)}({\mathbb R}^3) $,
that is, there is a constant $C>0$ such that 
\[
\| F(\psi _1)-F(\psi _2)\|_{H_{(s)}({\mathbb R}^3)} \leq C \|  \psi _1 - \psi _2 \|_{H_{(s)}({\mathbb R}^3)}
\left(\|  \psi _1   \|_{H_{(s)}({\mathbb R}^3)}^\alpha   +\|    \psi _2 \|_{H_{(s)}({\mathbb R}^3)}^\alpha  \right)
\]
{ for all $ \psi _1,\psi _2 \in H_{(s)}({\mathbb R}^3)$}.}

The polynomial in $\psi$  vector-valued functions and functions $F(\psi ) =(\gamma ^0 \psi, \psi) \gamma ^0 \psi $, $F(\psi ) =|\gamma ^0\gamma ^5 \psi|^\alpha   \psi $, $F(\psi ) =|\gamma ^0\gamma ^5 \psi|^\alpha  \gamma ^0 \psi $, $F(\psi ) =|\gamma ^0\gamma ^5 \psi|^\alpha  \gamma ^5 \psi $, $F(\psi ) = \pm|\psi |^{ \alpha } \psi $,  
are important examples of the Lipschitz continuous with exponent $\alpha  > 0$ in
the Sobolev space $H_{(s)}({\mathbb R}^3)$, { $s > 3/2$}, functions.

 We  consider the  Cauchy problem with the initial time $t=1 $. Furthermore,  we admit  that the mass of the field (particle)  can be changing in time and  vanishing at infinity. We allow the complex-valued mass $m$ keeping in mind the Huygens' principle that holds for the imaginary mass { \cite{ArX2021_HP}.}  (For the de~Sitter space see   \cite{JPh_2021}). More exactly, the model that we choose is determined by the Dirac operator
\begin{equation}
\label{DO}
{\mathscr{D}}_{FLRW}(t,\partial _t,\partial _x)
 := 
  \dsp 
 i {\gamma }^0    \partial_0   +i \frac{1}{a(t)}\sum_{j=1,2,3} {\gamma }^j  \partial_{j} +i \frac{3\dot a(t)}{2a(t)}    {\gamma }^0     -m t^{-1}{\mathbb I}_4  \,,
\end{equation}
where $m \in {\mathbb C}$.   This model includes  equation of neutrino with $m=0$ as well as massive spin-$\frac{1}{2}$ particle  undergoing a redshifting of its wavelength as the universe expands. { (See, e.g., \cite{Fock,Huang_Justin}.)}

Define the space
\[
X(R, k, \ell,m ) :=\Big\{\psi  (x,t)\in C([0,\infty);H_{(k)})\,\Big|\, \| \psi  \|_X :=  \sup_{t \in [1,\infty) }t^{3\ell/2-|\Im (m)|} \| \psi (t) \|_{H_{(k)}} \leq R   \Big\}\,.
\]
For the potential we  write { $V \in  {\mathscr{B}}^{(\ell,k)}$} that implies that all entries of the matrix $V$  belong to   the  space
\begin{eqnarray*} 
{\mathscr{B}}^{(\ell,k)}:= \Big\{v \in C_{t,x}^{\ell,k}  ([0,\infty)\times {\mathbb R}^3)  \,\Big|\, 
 \partial_t^j \partial_x^\alpha  v (x,t) \in L^\infty([0,\infty)\times {\mathbb R}^3)  , \,\forall \alpha ,\,|\alpha |\leq k,\,\,  
\forall j \leq \ell  \Big\}\,.
\end{eqnarray*}In the next theorem for the simplicity, we set $a_0=1$, that is,  $a(t)=t^{\ell}$.

\begin{theorem}
\label{TTFLRW_1}
Let $F=F(\psi ) \in C^3({\mathbb C}^4;{\mathbb C}^4)$  be the Lipschitz continuous with exponent $\alpha >0 $ in the space $H_{(k)}({\mathbb R}^3) $, $k\geq  3 $, function and  the potential $V \in {\mathscr{B}}^{(0,k)}$ is  self-adjoint, 
$V (x,t)= V ^*  (x,t) $. The problem
\begin{equation}
\label{0.6}
\cases{ \dsp 
\left(     \partial_0   +   t^{-\ell}\sum_{j=1,2,3}{\alpha }^j \partial_j+      \frac{3\ell}{2}t^{-1}{\mathbb I}_4     +i mt^{-1} {\gamma }^0 -iV (x,t)\right)\psi=F(\psi ) \,,\cr 
\psi (x,1)=\psi_0 (x )\,,
}
\end{equation} 
with 
\begin{equation}
\label{conell}
3\ell >2|\Im (m)|  \quad and  \quad  \alpha > 2\left(3\ell -2|\Im (m)| \right)^{-1}
\end{equation}
for sufficiently small $\varepsilon $ and $\psi_0 \in H_{(k)} $, {$k\geq  3  $}, $  \|\psi_0 \|_{(k)}\leq \varepsilon $, has a global solution $\psi \in X(2\varepsilon , k, \ell, m ) $. 
\end{theorem}

 For the real-valued  mass,   $m \in \R $, the  condition is $\alpha >2/(3\ell)$. In particular, 
$\alpha > 1$  if $\ell=2/3$,  while 
$\alpha > \frac{4}{3} $  if $\ell=1/2$. 
It is interesting, that for the Klein-Gordon  equation  with $m=0$ for $\ell=2/3$  if the semilinear term is  positive and the Lipschitz continuous with exponent $\alpha <2 $, then, the blowup occurs (see, \cite{Galstian-Yagdjian2020}).    The same is true, if $\ell =1/2$ and $m=0$ if $\alpha \leq 4/3$. For the general $\ell$ the blowup occurs if $\alpha <2/(3-3\ell)$ or $ \alpha <(4-\ell)/(4-3\ell)$. Later in \cite{Palmieri20, Palmieri21, Tsutaya} various proofs of this result were proposed.
\medskip

The equation (\ref{0.6}) is a symmetric hyperbolic system, and the local existence of the solution is known (see, e.g., \cite{Taylor_III}). Thus, the local Cauchy problem for (\ref{0.6}) is well posed  in $C^0([1,T];$ $(H_s({\mathbb R}^3))^4)  $, $s\geq 3$,  for some $T>1$ (see, e.g., \cite{Kato}). 
The next theorem states that some local solutions of the large data can be continued  to the global ones.

\begin{theorem}
\label{T4.1}
Let $ m \in \R$,   the potential $V \in {\mathscr{B}}^{(\infty,\infty)}$ is self-adjoint,  $V^*(x,t)=V(x,t)$, and  $V^T(x,t)\gamma ^2+\gamma ^2V(x,t)=0$. Suppose that the  scale factor    $a(t)=a_0t^{\ell}  $ has    
\begin{equation}
\label{aintehgral}
\ell >1/3 \,.
\end{equation}
Assume  also that  the function  $F=F(\xi ,\eta ) $, $F \in C^\infty({\mathbb R}^2;{\mathbb C}^4)$, has the form    
\begin{equation}
\label{F}
F\left(\psi^* \gamma ^0 \psi\,,\,\psi^* \gamma ^0\gamma ^5 \psi \right)= \alpha \left(\psi^* \gamma ^0 \psi\,,\,\psi^* \gamma ^0\gamma ^5 \psi \right)I +i\beta \left(\psi^* \gamma ^0 \psi\,,\,\psi^* \gamma ^0\gamma ^5 \psi \right) \gamma ^5,\,\, 
\end{equation}
where $\alpha $ and $\beta $ are real-valued functions,   with 
\begin{equation}
\label{Falphabeta} 
\alpha (\xi ,\eta ) =O(|\xi| +|\eta| ), \quad \beta  (\xi ,\eta ) =O(|\xi| +|\eta| ), \quad |\xi| +|\eta| \to 0\,. 
\end{equation}
 Assume also that the function $\Psi_0 =\Psi_0 (x)\in C_0^\infty ({\mathbb R}^3;{\mathbb C}^4)$ satisfies the Lochak-Majorana condition
\begin{equation}
\label{LMC}
\rho ^2(\Psi_0(x)):= | \Psi_0^*  (x) \gamma ^0\Psi_0 (x)|^2  + |\Psi_0^* (x) \gamma ^0\gamma ^5\Psi_0(x)|^2  = 0
\quad \mbox{\it for all }\quad x \in \R^3\,.  
\end{equation}

Then for $\chi _0 \in C_0^\infty ({\mathbb R}^3;{\mathbb C}^4) $
there is $\varepsilon _0>0 $ such that the Cauchy problem 
\begin{equation}
\label{NDE_CP}  
\cases{ \dsp\left( {\mathscr{D}}_{FLRW}(t,\partial _t,\partial _x) + \gamma ^0V(x,t)  \right)\psi    =F\left( \psi^*  \gamma ^0 \psi\,,\,\psi^* \gamma ^0\gamma ^5 \psi \right)\psi , \quad t>1  ,\cr
\psi (x,1)=\Psi _0 (x)+\varepsilon \chi _0 (x),
}
\end{equation}
 with  $ 0<\varepsilon <\varepsilon _0$ has a unique solution $\psi =\psi (x,t) $ such that { $\psi   \in  C^{1,\infty}_{t,x}([1,\infty) \times {\mathbb R}^3)   $}. The solution scatters to free solution  as $t \to +\infty$.
\end{theorem}

The last assertion of the previous theorem  follows from    Theorem~\ref{T4.1A} that states that if  
$F( \psi^* \gamma ^0 \psi,\psi^*\gamma ^0\gamma ^5 \psi )\psi$  is the Lipschitz continuous function with exponent $\alpha >0 $ in the space $H_{(3)}({\mathbb R}^3) $, 
 where $3\ell\alpha>2$, then for every solution $\psi =\psi (x,t) $ of the problem  
 (\ref{NDE_CP})
given by Theorem~\ref{T4.1} there exists 
a solution $\widetilde \psi (x,t)$ of  the free Dirac equation 
\[ 
\dsp\left( {\mathscr D}_{FLRW}(t,\partial_t,\partial_x)  + {\gamma }^0 V(x,t) \right)\widetilde \psi   =0  \,,
\]
such that 
\[
\lim_{t \to + \infty}   \left( \left\|   \psi (x,t)-   \widetilde \psi (x,t) \right\|_{(H_{(3)}({\mathbb R}^3))^4}+ \left\| \partial_t  \psi (x,t)-   \partial_t \widetilde \psi (x,t) \right\|_{(H_{(3)}({\mathbb R}^3))^4} \right)=0\,.
\]  
Moreover, if $V(x,t)=0 $, then
\[
\widetilde \psi  (x,t) 
 = {\mathscr{D}}^{co} (t,\partial _t,\partial _x)\gamma^0\left (
   \begin{array}{ccccc}
  {\cal K}_1(x,t,D_x;m;1 ) {\mathbb I}_2& {\mathbb O}_2   \\
   {\mathbb O}_2&  {\cal K}_1(x,t,D_x;-m;1 ){\mathbb I}_2   \\ 
   \end{array}
   \right) [\widetilde \psi (x,1)] ,\quad t>1,
\]
with the operators ${\mathscr{D}}^{co} (t,\partial _t,\partial _x) $ and ${\cal K}_1(x,t,D_x;m;1 ) $ written in the explicit way (see \cite{JDE2021} or subsection~\ref{SS7.3}).

The rest of this paper is organized as follows.  In Section~\ref{S1}, we prove Theorem~\ref{TTFLRW_1} by  the energy estimate (subsection~\ref{SSEE})  and fixed point argument (subsection~\ref{SS1.2}).
 The asymptotics on the positive half-line of time (Theorem~\ref{T1.2}) and the representation of the solution of the free  Dirac equation in the FLRW spacetime  are given in  subsection~\ref{SS7.3}. In subsection~\ref{SS5.4}, to make the text self-contained,  we prove that the classical solutions of the semilinear Dirac equation obey the finite propagation speed property. In Section~\ref{S2}, we analyze  the Lochak-Majorana condition in the FLRW spacetime and obtain its time-evolution. In Section~\ref{S3}, we prove   Theorem~\ref{T4.1} except the asymptotics part that follows from Theorem~\ref{T4.1A} of  
 Section~\ref{S4}. There, we also give the asymptotics at infinity for the solutions of Theorem~\ref{T4.1}. 
In  Section~\ref{S5}, we demonstrate  the blow-up result (Theorem~\ref{TBUH+}) for the solution of the Dirac equation in FLRW spacetime.

\section{Small data global existence}

\label{S1}

\subsection{Energy estimate}
\label{SSEE}
\setcounter{equation}{0}

The metric tensor  (\ref{metric})   allows us  to develop the energy estimate that is the main tool to prove global existence of small amplitude solutions. 
\begin{lemma}
\label{L1.2FLRW}
Assume that $a(t)=a_0t^\ell$ and $f \in C([0,\infty);(H_{(k)}({\mathbb R}^3))^4 )$, the potential $V \in {\mathscr{B}}^{(0,k)}$ is self-adjoint,  $V (x,t)= V^* (x,t) $.  
Then, for the solution of 
\[
  \dsp 
\left(     \partial_0   +   \frac{1}{a(t)}\sum_{j=1,2,3}{\alpha }^j  \partial_j+      \frac{3\dot a(t)}{2a(t)}{\mathbb I}_4     +i mt^{-1} {\gamma }^0 -iV (x,t)\right)\psi=f(x,t) \,,
\]
in the Sobolev space $H_{(k)}({\mathbb R}^3)$ one has   
\begin{eqnarray}
\label{ForwFLRW}
&  &
\|\psi (t) \|_k\\
& \leq  & 
c  t^{-3\ell/2+|\Im (m)|}s^{ 3\ell/2-|\Im (m)|}\|\psi (s) \|_k +
 c   t^{-3\ell/2+|\Im (m)|}\int_s^t \tau ^{ 3\ell/2-|\Im (m)|} \| f (\tau )\|_k \,d\tau   \,,\quad 1\leq s\leq t\,, \nonumber \\
\label{BackFLRW}
&  &
\|\psi (t ) \|_k \\
& \leq  &
c  t^{-3\ell/2-|\Im (m)|}s^{ 3\ell/2+|\Im (m)|}\|\psi (s ) \|_k  +  c t^{-3\ell/2-|\Im (m)|}\int_t^s  \tau ^{ 3\ell/2+|\Im (m)|}\| f (\tau )\|_k \, d\tau \,,\quad 1\leq  t\leq s .  \nonumber
\end{eqnarray}
\end{lemma}
\medskip

\ndt
{\bf Proof.} Consider the energy integral
\begin{eqnarray*}
E(t)
& = &
\int_{{\mathbb R}^3}|\psi (x,t)|^2\,dx \,.
\end{eqnarray*} 
Then, 
\[
\frac{d}{d t} E(t)
  =  \int_{{\mathbb R}^3}\left(  -  3\ell t^{-1} |\psi (x,t)|^2  + 2t^{-1}(\Im (m)) \psi^* (x,t)  {\gamma }^0 \psi (x,t)  \right) \,dx
 +\int_{{\mathbb R}^3}  2\Re(  f^* (x,t)  \psi (x,t)  ) \,dx\,.
\]
Furthermore, 
\begin{eqnarray*}
&  &
  - 2t^{-1}|\Im (m)|  |\psi (x,t)|^2 \leq 
  2t^{-1}(\Im (m))\psi^* (x,t)  {\gamma }^0 \psi (x,t) \leq 
 2t^{-1}|\Im (m)|  |\psi (x,t)|^2 \,.
\end{eqnarray*}
Denote
\[
\delta _+(t):=    (- 3\ell    + 2|\Im (m)|)t^{-1},\qquad \delta _-(t):=    (- 3\ell    - 2|\Im (m)|)t^{-1}\,.
\]
Then
\[
\delta _-(t)   |\psi (x,t)|^2 \leq 
- 3\ell t^{-1} |\psi (x,t)|^2  + 2t^{-1}(\Im (m))\psi^* (x,t)  {\gamma }^0 \psi (x,t)  \leq 
\delta _+ (t) |\psi (x,t)|^2 
\]
and
\begin{equation}
\label{DINFLRW}
\delta _- (t)  E(t)
 +\int_{{\mathbb R}^3}2\Re(  f^* (x,t)  \psi (x,t)  ) dx 
\leq  
\frac{d}{d t} E(t)\leq 
\delta _+ (t)  E(t)
 +\int_{{\mathbb R}^3} 2\Re(  f^* (x,t)  \psi (x,t)  ) dx  .
\end{equation}
In particular, 
\begin{eqnarray*}
\frac{d}{d t} \left( E(t)\exp \left( \dsp -\int_1^t \delta _+(\tau )\,d\tau \right)  \right)
& \leq  & 
 2\exp \left( \dsp -\int_1^t \delta _+(\tau )\,d\tau \right)\| f (x,t)\| \| \psi (x,t) \|\,.
\end{eqnarray*}
We integrate the last inequality from $s$ to $t$, $s<t$, and obtain
\begin{eqnarray}
\label{EIN}
 E(t)\exp  \left( \dsp -\int_1^t \delta _+(\tau )\,d\tau \right)  
& \leq  &
 E(s)\exp  \left(\dsp -\int_1^s \delta _+(\tau )\,d\tau \right) \\
&  &
+  2\int_s^t \,\exp  \left(\dsp -\int_1^\tau  \delta _+(\tau_1 )\,d\tau_1 \right) \| f (x,\tau )\| \| \psi (x,\tau ) \|\,d \tau . \nonumber
\end{eqnarray}
We fix $s$ and for $k=0$   denote
\[
y(t) :=\max_{\tau \in [s,t]}\exp \left( \dsp   -\frac{1}{2}\int_1^\tau  \delta _+(\tau_1 )\,d\tau_1 \right)  \|\psi (\tau ) \|_k, 
\]
then
\[
  y^2(t) =\max_{\tau \in [s,t]} \exp \left( \dsp   -\int_1^\tau  \delta _+(\tau_1 )\,d\tau_1 \right)\|\psi (\tau ) \|_k^2.  
\]
Hence (\ref{EIN}) implies
\[
y^2(t)
 \leq  
y^2(s)+
 2y(t)\int_s^t \exp \left(\dsp   -\frac{1}{2}\int_1^\tau  \delta _+(\tau_1 )\,d\tau_1\right) \| f (\tau )\| \,d\tau \,.
\]
Furthermore,
\[
y (t)
\leq  
c  y (s)+
 c     \int_s^t  \exp \left(\dsp   -\frac{1}{2}\int_1^\tau  \delta _+(\tau_1 )\,d\tau_1\right)\| f (\tau )\| \,d\tau  \,.
\]
Thus, we   have 
\[
\|\psi (t) \|_k 
\leq  
c  \exp \left(\dsp   \frac{1}{2}\int_s^t  \delta _+(\tau_1 )\,d\tau_1\right) \|\psi (s) \|_k +
 c   \int_s^t  \exp \left(\dsp  \frac{1}{2}\int_ \tau ^t \delta _+(\tau_1 )\,d\tau_1 \right) \| f (\tau )\| \,d\tau   \,, 
\]
for all {$1\leq s\leq t $}. Next we choose the left hand side of (\ref{DINFLRW}), rewrite it as
\begin{eqnarray*}
\exp \left(\dsp -\int_1^t \delta _-(\tau )\,d\tau \right) \int_{{\mathbb R}^3} 2\Re(  f^* (x,t)  \psi (x,t)  ) \,dx
& \leq  &
\frac{d}{d t} \left( E(t)\exp \left(\dsp -\int_1^t \delta _-(\tau )\,d\tau \right)\right) \,, 
\end{eqnarray*}
integrate from $t$ to $s$, $s\geq t\geq 1 $, and obtain 
\begin{eqnarray*}
&  &
- 2\int_t^s \exp \left(\dsp -\int_1^\tau  \delta _-(\tau_1 )\,d\tau_1 \right) \| f (\tau )\| \| \psi (\tau ) \|\, d\tau \\
 & \leq &  
  E(s)\exp \left(\dsp -\int_1^s  \delta _-(\tau_1 )\,d\tau_1 \right)  
-   E(t)\exp \left(\dsp -\int_1^t  \delta _-(\tau_1 )\,d\tau_1 \right)  
\end{eqnarray*}
that implies
\begin{eqnarray*}
&  &
E(t)\exp \left(\dsp -\int_1^t  \delta _-(\tau_1 )\,d\tau_1 \right)  \\
& \leq  &
  E(s)\exp \left(\dsp -\int_1^s  \delta _-(\tau_1 )\,d\tau_1 \right) +  2\int_t^s \exp \left(\dsp -\int_1^\tau  \delta _-(\tau_1 )\,d\tau_1\right)\| f (\tau )\| \| \psi (\tau ) \|\, d\tau  \,,\quad 1\leq  t\leq s \,.
\end{eqnarray*}
If we  fix $s$ and for $k=0$ denote
\[
z(t) :=\max_{\tau \in [t,s]} \exp \left(\dsp - \frac{1}{2}\int_1^\tau   \delta _-(\tau_1 )\,d\tau_1\right)\|\psi (\tau ) \|_k, 
\]
then
\[  z^2(t) =\max_{\tau \in [t,s]} \exp \left(\dsp - \int_1^\tau   \delta _-(\tau_1 )\,d\tau_1 \right)\|\psi (\tau ) \|_k^2\,.  
\]
The inequality
\begin{eqnarray*}
 z^2(t) 
& \leq  &
z^2(s) + 2 z(t) \int_t^s \exp \left(\dsp -\frac{1}{2}\int_1^\tau   \delta _-(\tau_1 )\,d\tau_1 \right)\| f (\tau )\|_k \, d\tau \,,\quad 1\leq t\leq s, 
\end{eqnarray*}
implies
\begin{eqnarray*}
z (t) 
& \leq  &
cz(s)  +  c \int_t^s \exp \left(\dsp -\frac{1}{2}\int_1^\tau   \delta _-(\tau_1 )\,d\tau_1 \right)\| f (\tau )\|_k \, d\tau \,,\quad 1\leq t\leq s, 
\end{eqnarray*}
and, consequently, 
\begin{eqnarray*}
&  &
\exp \left(\dsp -\frac{1}{2}\int_1^t   \delta _-(\tau_1 )\,d\tau_1 \right)\|\psi (t ) \|_k\\
& \leq  &
c\exp \left(\dsp -\frac{1}{2}\int_1^s   \delta _-(\tau_1 )\,d\tau_1 \right)\|\psi (s ) \|_k  +  c \int_t^s \exp \left(\dsp -\frac{1}{2}\int_1^\tau   \delta _-(\tau_1 )\,d\tau_1 \right)\| f (\tau )\|_k \, d\tau \,,
\end{eqnarray*}
for all $1\leq t\leq s $. Thus, for all $1\leq  t\leq s$ we have 
\begin{eqnarray*}
\|\psi (t ) \|_k
& \leq  &
c\exp \left(\dsp  \frac{1}{2}\int_s^t   \delta _-(\tau_1 )\,d\tau_1 \right)\|\psi (s ) \|_k  +  c \int_t^s \exp \left(\dsp \frac{1}{2}\int_\tau^t    \delta _-(\tau_1 )\,d\tau_1  \right)\| f (\tau )\|_k \, d\tau \,.
\end{eqnarray*}
For every $k > 0$ the inequalities of the lemma can be obtained  from the case of $k = 0$ by differentiation.
The lemma is proved.  \qed

\subsection{Global existence. Small data. Proof of Theorem~\ref{TTFLRW_1}}
\label{SS1.2}

Set $s=1$ in (\ref{ForwFLRW}), then due to the condition (${\mathcal L}$) for $k>3/2$   and all $t >1$  one has
\begin{equation}
\label{1.5b}
\|\psi (t) \|_k
\leq  
c  t^{-3\ell/2+|\Im (m)|}\|\psi (1) \|_k +
 c  t^{-3\ell/2+|\Im (m)|} \int_1^t  \tau  ^{3\ell/2-|\Im (m)|}\| F(\psi ) (\tau )\|_k  \,d\tau  . 
\end{equation}
Hence, 
\[
\sup_{t \in [1,\infty) } t^{3\ell/2-|\Im (m)|} \| \psi (t) \|_{k} \\
 \leq  
 c\|\psi (1) \|_{(s)}+  c\left( \sup_{t \in [1,\infty) }t^{3\ell/2-|\Im (m)|} \| \psi (t) \|_{k} \right)^{1+\alpha }\int_1^t \left(\tau ^{3\ell/2-|\Im (m)|}\right)^{-\alpha }   \,d\tau .
\]
The last integral  is
\begin{eqnarray*}
\int_1^t \tau ^{-\alpha (3\ell/2-|\Im (m)|)}  \,d\tau 
& = &
\frac{1}{1-\alpha(\frac{3}{2}\ell -|\Im (m)|)}\left( t^{1-\alpha(\frac{3}{2}\ell -|\Im (m)|) }-1 \right)\,.
\end{eqnarray*}
Hence, the integral is bounded for all $t \in [1,\infty)$ because of  the condition  (\ref{conell}). 

Denote by $S(t,s) $ the propagator (fundamental solution for the Cauchy problem),  that is, an operator-valued  solution of 
the problem 
\begin{eqnarray}
\label{Prop}
\hspace{-0.9cm} &  &
\cases{   \dsp 
\left(     \partial_0   +  t^{-\ell }\sum_{j=1,2,3}{\alpha }^j  \partial_j+      \frac{3\ell}{2}t^{-1}{\mathbb I}_4     +i mt^{-1} {\gamma }^0 -iV (x,t)\right)S(t,s)=0\,,\quad t , s \in \R_+\,,\cr 
S(s,s)=I \,(\mbox{identity operator})\,.
}
\end{eqnarray}
Then the solution of the problem  
\begin{eqnarray*}
&  &
\cases{  \dsp  
\left(     \partial_0   +  t^{-\ell }\sum_{j=1,2,3}{\alpha }^j  \partial_j+      \frac{3\ell}{2}t^{-1}{\mathbb I}_4     +i mt^{-1} {\gamma }^0 -iV (x,t)\right)\psi(x,t)=f(x,t)\,,\quad t \geq s >0\,,\cr
\psi (x,s)=\psi_0 (x )\,,
} 
\end{eqnarray*}
is given by Duhamel's principle
\[
\psi (x,t)= S(t,s) \psi_0 (x )+\int_s^t S(t,\tau )f(x,\tau)\,d \tau\,.
\]
It is known (see, e.g., \cite{Kato,Taylor_II}) that for $ \psi_0 \in H_{(k)}({\mathbb R}^3)$ and $f \in C([0,\infty)  ; H_{(k)}({\mathbb R}^3))$, $k>5/2$, the unique solution $\psi  \in C([1,\infty)  ; H_{(k)}({\mathbb R}^3))\cap C^1([1,\infty)  ; H_{(k-1)}({\mathbb R}^3))$ exists.  
Next, we define the operator $\mathcal S$ by
\[
{\mathcal  S} \psi(x,t) :=S(t,1) \psi _0 ( x)+\int_1^t S(t,\tau )F(\psi(x,\tau))\,d \tau\,.
\]
We  are going to prove that ${\mathcal  S}$ is a contraction   
\begin{eqnarray*}
&  &
{\mathcal  S} : X(R, k, \ell,m ) \longrightarrow  X(R, k, \ell,m )
\end{eqnarray*}
for sufficiently small $R$.
The  inequality (\ref{1.5b}) proves that the operator $S$ 
maps $X(R, k, \ell,m ) $ 
 into itself if $ \psi _0 \in H_{(k)}$, $\|\psi _0\|_{H_{(k)}}< \varepsilon $, and   $\varepsilon  $ and $R$ are sufficiently small, namely, if 
\[
\varepsilon +  \left(\frac{1}{\alpha(\frac{3}{2}\ell -|\Im (m)|)-1}
\right)  R^{1+ \alpha }  < R .
\]

In order to verify a contraction property, we write
\begin{eqnarray*}
&  &
{\mathcal  S} \psi_1(x,t) - {\mathcal  S} \psi_2(x,t) = \int_1^t S(t,\tau )\left(  F(\psi_1(x,\tau)) -F(\psi_2(x,\tau)) \right) \,d \tau
\end{eqnarray*} 
and use the condition ($\mathcal L$) and Lemma~\ref{L1.2FLRW} to estimate the norm 
\begin{eqnarray*} 
&  &
\|{\mathcal  S} \psi_1(x,t) - {\mathcal  S} \psi_2(x,t)\|_k  \\
& \leq  &
 c   t^{-3\ell/2+|\Im (m)|}  \int_1^t  \tau ^{3\ell/2-|\Im (m)|}\|   F(\psi_1( \tau)) -F(\psi_2( \tau)) \|_k \,d\tau  \\ 
& \leq  & 
 c t^{-3\ell/2+|\Im (m)|} \int_1^t \tau ^{3\ell/2-|\Im (m)|} 
\|  \psi _1 (\tau)- \psi _2(\tau) \|_k
\left(\|  \psi _1 (\tau)  \|_k^\alpha   
+\|    \psi _2 (\tau)\|_k^\alpha  \right)\,d \tau \,.
\end{eqnarray*} 
It follows
\begin{eqnarray*} 
&   & 
\sup_{t \in [1,\infty) }  t^{3\ell/2-|\Im (m)|}\|{\mathcal  S} \psi_1(x,t) -{\mathcal  S} \psi_2(x,t)\|_k \\ 
& \leq  &
c   \left( \sup_{\tau  \in [1,\infty) } \tau ^{3\ell/2-|\Im (m)|}
\|  \psi _1 (\tau)- \psi _2(\tau) \|_k \right)  \int_1^t 
\left(\|  \psi _1 (\tau)  \|_k^\alpha   
+\|    \psi _2 (\tau)\|_k^\alpha  \right)\,d \tau \\ 
& \leq  &
c   \left( \sup_{\tau  \in [1,\infty) } \tau ^{3\ell/2-|\Im (m)|}
\|  \psi _1 (\tau)- \psi _2(\tau) \|_k \right) \\
&  &
\times   \left( \sup_{\tau  \in [1,\infty) } \tau ^{3\ell/2-|\Im (m)|}
\left(\|  \psi _1 (\tau)  \|_k   
+\|    \psi _2 (\tau)\|_k  \right)\right)^\alpha   \int_1^t 
\tau ^{-\alpha (3\ell/2-|\Im (m)|)}\,d \tau \\  
& \leq  &
 C {  \|\psi _1 -\psi _2 \|_{X(R, k, \ell,m )}}
R^\alpha \,.
\end{eqnarray*} 
Here we have used condition (\ref{conell}). Then we choose  $R$ such that $CR^{\alpha} <1 $. The Banach fixed-point theorem completes the proof of Theorem~\ref{TTFLRW_1}.  \qed

\subsection{Large  time asymptotics }   
\label{SS7.3}

For the equation with    $\ell\not= 1$  we define 
the forward and backward light cones   as the boundaries of  $D_+(x_0, t_0)$ and $D_-(x_0, t_0)$, respectively, where
\begin{equation}
\label{Dpm} 
D_{\pm}\left(x_{0}, t_{0}\right) :=\left\{(x, t) \in {\mathbb R}^{3+1} ;
\left|x-x_{0}\right| \leq \pm\left(\phi (t) -\phi (t_{0})  \right)\right\}\,,
\end{equation}
and  $\phi (t):= \frac{1}{1-\ell}t^{1-\ell}$. 
Any intersection of $D_-(x_0, t_0)$ with the hyperplane $ t = const  $, $0<  const < t_0$, determines the
so-called {\it dependence domain} for the point $ (x_0, t_0)$, while the intersection of $D_+(x_0, t_0)$
with the hyperplane $t = const$, $ const> t_0>0$ is the so-called {\it domain of influence} of the point
$ (x_0, t_0)$. 
The Dirac equation (\ref{DE})  is non-invariant with respect to time inversion and its solutions have different properties in different directions of time.

A retarded fundamental solution (a retarded inverse) for the Dirac operator (\ref{DO})  is a four dimensional matrix with the 
distribution-valued  entries $ 
{\mathcal E}^{ret}={\mathcal E}^{ret} \left(x, t ; x_{0}, t_{0};m\right)
 $      that solves the equation
\begin{equation}
\label{FSE} 
{\mathscr{D}}_{FLRW}(t,\partial _t,\partial _x){\mathcal E}  \left(x, t ; x_{0}, t_{0};m\right) 
 = 
\delta ( x-x_0) \delta (t-t_0) {\mathbb I}_4, \qquad 
 (x,t), (x_0,t_0 ) \in {\mathbb R}^3\times {\mathbb R}_+,  
\end{equation}
and with the support in the {\it chronological future} (causal future) $D_+(x_0, t_0)$  of the point $(x_0,t_0)  \in {\mathbb R}^3\times {\mathbb R}_+$. The
advanced fundamental solution (propagator) $ 
{\mathcal E}^{adv}={\mathcal E}^{adv}  (x, t ; x_{0}, t_{0};$ $  m )
 $ solves the equation (\ref{FSE}) and has the  support in the {\it chronological past} (causal
past) $D_-(x_0, t_0)$.  
\medskip

Following \cite{JDE2021}    we introduce two more $\gamma$-matrices (projection operators),
 the upper left corner and lower right corner matrices,
\begin{eqnarray*}
\gamma^U
& = &
\left(
\begin{array}{cc}
 {\mathbb I}_2 & {\mathbb O}_2   \\
 {\mathbb O}_2 & {\mathbb O}_2   \\ 
\end{array}
\right)=\frac{1}{2}({\mathbb I}_4+ \gamma^0 )
,\qquad 
 \gamma^L
  =  \left(
\begin{array}{cc}
 {\mathbb O}_2 & {\mathbb O}_2   \\
 {\mathbb O}_2 & {\mathbb I}_2   \\ 
\end{array}
\right)=
\frac{1}{2}({\mathbb I}_4- \gamma^0 )\,.
\end{eqnarray*}
Next we define the right co-factor
\begin{eqnarray*}
{\mathscr{D}}^{co}_{FLRW} (t,\partial _t,\partial _x)
& := &
i t^{-\frac{\ell }{2}} \gamma ^0 (t^{ i m}  \gamma^U+t^{-i m}\gamma^L )\frac{\partial }{\partial t} 
+i t^{-\frac{3\ell }{2}}  \sum_{k=1}^3\gamma ^k (t^{ i m}  \gamma^U+t^{-i m}\gamma^L )\frac{\partial }{\partial x_k} 
\end{eqnarray*}
 of the Dirac operator ${\mathscr{D}}_{FLRW} (t,\partial _t,\partial _x) $ of (\ref{DO}). The composition   
${\mathscr{D}}_{FLRW} (t,\partial _t,\partial _x)$ ${\mathscr{D}}^{co}_{FLRW} (t,\partial _t,\partial _x) $ is a diagonal matrix of operators (see \cite{JDE2021}).

Denote  ${\mathcal E}^w(x,t)$  the fundamental solution to the Cauchy problem for the  wave equation  in the Minkowski space, that is, the solution of the problem
\[
{\mathcal E}^w_{ tt} -   \Delta  {\mathcal E}^w  =  0 \,, \quad {\mathcal E}^w(x,0)=\delta (x)\,, \quad {\mathcal E}^w_{t}(x,0)= 0\,.
\]
Here $\Delta  $ is the Laplace operator in ${\mathbb R}^3$. 
Henceforth, $F \left(  \alpha , \beta ;\gamma ;z\right) $ is the hypergeometric function (see, e.g., \cite{B-E}). 
For $n=3$ 
  (see, e.g., \cite{Shatah})  
\[ 
{\mathcal E}^{w}(x, t) :=\frac{1}{4\pi } \frac{\partial}{\partial t}
\frac{1}{t} \delta(|x|-t)\,.
\]
The distribution $\delta(|x|-t)$ is defined
 by
\[
\langle \delta(|\cdot|-t), \psi(\cdot)\rangle  =\int_{|x|=t} \psi(x)\, d x \quad \mbox{\rm for } \quad \psi \in C_{0}^{\infty}\left({\mathbb R}^{3}\right).
\]

According to Theorem~1.1~\cite{JDE2021},
for every    $x_0 \in {\mathbb R}^3$, $ t, t_0 \in {\mathbb R}_+$, the retarded fundamental solution  of the Dirac operator   ${\mathcal E}^{ret}(x ,t;x_0 ,t_0;m) $ 
is given  as follows
\begin{eqnarray*}
&  &
{\mathcal E}^{ret}(x ,t;x_0 ,t_0;m) \\
& = &
- 2  t_0^{\frac{\ell }{2}-i m } {\mathscr{D}}^{co} (t,\partial _t,\partial _x)   \int_0^{ \phi (t)- \phi (t_0 ) }\left (
   \begin{array}{ccccc}     
 E(r,t;t_0; m)    {\mathbb I}_2 & {\mathbb O}_2 \\
   {\mathbb O}_2&   
 E(r,t;t_0; -m)   {\mathbb I}_2   \\ 
   \end{array}
   \right) {\mathcal E}^w(x-x_0,r)\,dr\,.
\end{eqnarray*}
Here the kernel is defined as follows (see \cite{JDE2015,JDE2021})
\begin{eqnarray*}
\hspace{-0.3cm}  E(r,t;t_0  ; m) 
& := & 
{ 2^{ \frac{2 im}{ 1-\ell }-1}} (1-\ell )^{\frac{\ell }{1-\ell }}   \phi (t_0  )^{\frac{ \ell+2im}{1-\ell }}
\left( \left(\phi (t)+ \phi (t_0  )  \right)^2- r^2 \right)^{-  \frac{ im}{ 1-\ell }} \\
\hspace{-0.3cm} &  &
\times  F \left(  i \frac{m}{ 1-\ell },  i\frac{m}{ 1-\ell };1;\frac{\left( \phi (t)-  \phi (t_0  )   \right)^2-r^2}{\left( \phi (t)+ \phi (t_0  ) \right)^2-r^2}\right)  \,.  \nonumber 
\end{eqnarray*}
The fundamental solution ${\mathcal E}_{+}(x ,t;x_0 ;m;\varepsilon )$ to the Cauchy problem for the Dirac equation is given by Theorem~1.2 \cite{JDE2021}. It states that  
${\mathcal E}_{+}(x ,t;x_0 ;m;\varepsilon )$, that is, a distribution satisfying 
\begin{eqnarray*}
&  &
\cases{{\mathscr{D}}_{FLRW}(t,\partial _t,\partial _x)
{\mathcal E}_{+}(x ,t;x_0 ;m;\varepsilon )=   {\mathbb O}_4 ,  \cr
{\mathcal E}_{+}(x ,\varepsilon ;x_0 ;m;\varepsilon )=\delta ( x-x_0){\mathbb I}_4\,, 
}
\end{eqnarray*}
is given by
\begin{eqnarray*}
{\mathcal E}_{+}(x ,t;x_0 ;m;\varepsilon )
& = &
-i\varepsilon ^{1+\frac{\ell }{2} -i m} (1-\ell)^{-1}{\mathscr{D}}^{co}_{FLRW} (x,t,\partial _t,\partial _x)\gamma^0\\
&  &
\times \int_0^{\phi (t)- \phi (\varepsilon )} \left (
   \begin{array}{ccccc}
K_1 \left(r,t; m ;\varepsilon \right)  {\mathbb I}_2& {\mathbb O}_2   \\
   {\mathbb O}_2& K_1 \left(r,t; -m ;\varepsilon \right) {\mathbb I}_2   \\ 
   \end{array}
   \right){\mathcal E}^w(x-x_0,r) \,dr     ,
\end{eqnarray*} 
where
\[
 K_1 \left(r,t; m ;\varepsilon \right) :=
 2^{2 i  \frac{m}{1-\ell}} \phi (\varepsilon )^{2 i  \frac{m}{1-\ell}-1}   \left( \left(\phi (t)+ \phi (\varepsilon )  \right)^2-  r ^2\right)^{-  i \frac{m}{1-\ell}}   F \left(  i \frac{m}{1-\ell},  i \frac{m}{1-\ell};1;\frac{\left(\phi (t)- \phi (\varepsilon )\right)^2- r ^2}{\left(\phi (t)+ \phi (\varepsilon )   \right)^2-  r ^2 }\right). \nonumber 
\] 
In order to write the solution to   the Cauchy problem  
  we introduce the operator  
\begin{eqnarray*}
&  &
{\cal G}(x,t,D_x;m )[f] (x,t)\\
& : =  &
- 2\int_\varepsilon ^{t  }b^{\frac{\ell }{2}-i m } \,d b\int_0^{ \phi (t)- \phi (b ) }  
 E(r,t;b ; m)\int_{{\mathbb R}^3} {\mathcal E}^w(x-y,r)  f( y,b)  \,dy \,dr
, \quad  f \in C_0^\infty({\mathbb R}^{n+1}),
\end{eqnarray*}
and  the operator ${\cal K}_1(x,t,D_x;m;\varepsilon)$   as follows:
\begin{eqnarray*}
&  &
{\cal K}_1(x,t,D_x;m;\varepsilon) [\varphi](x,t) \\
& := &
-i\varepsilon ^{1+\frac{\ell }{2} -i m}(1-\ell)^{-1}
\int_0^{ \phi (t)- \phi (\varepsilon ) }   K_1 \left(r,t; m ;\varepsilon \right)\int_{{\mathbb R}^3} {\mathcal E}^w(x-y,r)  \varphi (y )\,dy\,dr  
, \quad \varphi \in C_0^\infty({\mathbb R}^n).
\end{eqnarray*}
According to Theorem~1.3~\cite{JDE2021} the solution to the Cauchy problem
\begin{eqnarray*}
\cases{
{\mathscr{D}}_{FLRW}(t,\partial _t,\partial _x)\Psi  (x,t)=F( x,t)\,,\quad t>\varepsilon >0\,,\cr
\Psi (x,\varepsilon )=\Psi _\varepsilon (x)\,,
}
\end{eqnarray*}
with $m \in {\mathbb C}$, is given as follows
\begin{eqnarray*}
\Psi  (x,t)
& = &
 {\mathscr{D}}^{co}_{FLRW} (t,\partial _t,\partial _x)\Bigg\{\left (
   \begin{array}{ccccc}
    {\cal G}(x,t,D_x;m ) {\mathbb I}_2 & {\mathbb O}_2 \\
   {\mathbb O}_2&  {\cal G}(x,t,D_x;-m ){\mathbb I}_2   \\ 
   \end{array}
   \right) [F]  (x,t)\\
   &  &
 + \gamma^0\left (
   \begin{array}{ccccc}
  {\cal K}_1(x,t,D_x;m;\varepsilon ) {\mathbb I}_2& {\mathbb O}_2   \\
   {\mathbb O}_2&  {\cal K}_1(x,t,D_x;-m;\varepsilon ){\mathbb I}_2   \\ 
   \end{array}
   \right) [\Psi  _{\varepsilon }]  (x,t)  \Bigg\},\quad t> \varepsilon >0\,.
\end{eqnarray*}
Let $S(t,s) $ be the propagator (fundamental solution for the Cauchy problem) of the Dirac  equation, that is, an operator-valued  solution of 
the problem (\ref{Prop}).
Hence,
\[
\psi (x,t)=S(t,s)\psi (x,s)\,,\qquad t,s \in (0,\infty),\quad x\in \R^3\,.
\]
For the case of $V(x,t)=0$ the operator $S(t,s ) $  is written in the explicit form in \cite{JDE2021}.

\begin{theorem}
\label{T1.2}
Assume that
\[
2|\Im (m)|+\alpha \left(|\Im (m)| -\frac{3}{2}\ell \right)  <-1 \,.
\]
Let the function $\psi =\psi (x,t )  { \in X(2\varepsilon , k, \ell, m )}$, { $k\geq 3$}, be the solution of the problem (\ref{0.6}) given by Theorem~\ref{TTFLRW_1}.  
Then the  limit
\[
\lim_{t \to  \infty}   \int_1^{ t} S(1,\tau ) F(\psi (x,\tau ))\,d \tau   
\]
exists in the space  $(H_{(k)}({\mathbb R}^3))^4 $. Furthermore, the solution $\widetilde \psi (x,t)$ of the Cauchy problem for the free Dirac equation 
\begin{eqnarray*} 
\cases{ \dsp\left( {\mathscr D}_{FLRW}(t,\partial_t,\partial_x)  + {\gamma }^0 V(x,t) \right)\widetilde \psi (x,t)  =0  \,,\cr
\widetilde\psi (x,1)= \psi_0^+ (x)\,,
}
\end{eqnarray*}
where 
\[
 \psi_0^+ (x)={ \psi  _0 } (x)+\lim_{t \to  \infty}   \int_1^{ t} S(1,\tau ) F(\psi (x,\tau ))\,d \tau 
\]
satisfies 
\begin{eqnarray*}
\lim_{t \to  \infty}    \left\|  \psi (x,t)-   \widetilde \psi (x,t) \right\|_{(H_{(k)}({\mathbb R}^3))^4}=0 \,,
\end{eqnarray*} 
and { ${\mathscr S}\,:\, \psi \longmapsto \widetilde\psi  $ is a continuous  operator in $ {   X(2\varepsilon , k, \ell, m )}$, while  ${\mathscr S}_0\,:\, \psi_0 \longmapsto \psi ^+ $   is a continuous  operator in  $H_{(k)}({\mathbb R}^3) $.}

Moreover, if $V(x,t)=0 $, then
\[
\widetilde \psi  (x,t) \\
  =   {\mathscr{D}}^{co} (t,\partial _t,\partial _x)\gamma^0\left (
   \begin{array}{ccccc}
  {\cal K}_1(x,t,D_x;m;1 ) {\mathbb I}_2& {\mathbb O}_2   \\
   {\mathbb O}_2&  {\cal K}_1(x,t,D_x;-m;1 ){\mathbb I}_2   \\ 
   \end{array}
   \right) [ \psi_0^+ ]  (x,t),\quad t\geq 1\,.
\]  
\end{theorem}
\medskip

\noindent
{\bf Proof.} According to Lemma~\ref{L1.2FLRW} 
\begin{eqnarray*}
&  &
\|S(1,\tau )\psi (\tau  )  \|_k \leq 
c \tau ^{3\ell/2 +|\Im (m)|} 
\|\psi (\tau  ) \|_k ,  \quad  1 \leq \tau  \,.
\end{eqnarray*}
Then an existence of the limit follows from 
\begin{eqnarray*}
\| S(1,\tau ) F(\psi (x,\tau ))\|_k  
& \leq  &
c \tau ^{3\ell/2 +|\Im (m)|}\| F(\psi (x,\tau )) \|_k   \\
& \leq  &
c \tau ^{3\ell/2 +|\Im (m)|}\| \psi (x,\tau ) \|_k^ {1+\alpha}    \\
& \leq  &
c \tau ^{3\ell/2 +|\Im (m)|}\left[ \tau ^{-3\ell/2+|\Im (m)|}\right]^{1+\alpha }   \\
& \leq  &
c \tau ^{|\Im (m)|(2+\alpha)-3\ell \alpha /2} \,,
\end{eqnarray*}
where
\begin{eqnarray*}
&  &
|\Im (m)|(2+\alpha) -\frac{3}{2}\ell \alpha = 2|\Im (m)|+\alpha \left(|\Im (m)| -\frac{3}{2}\ell \right)  <-1
\end{eqnarray*}
according to the condition of the theorem. 

{ To prove continuity of the operator  ${\mathscr S}\,:\, \psi \longmapsto \widetilde\psi   $ consider two solutions $\psi_1= \psi_1(x,t)     $ and $\psi_2=\psi_2 (x,t) $ given by Theorem~\ref{TTFLRW_1} and the corresponding functions $\widetilde \psi_1 (x,t) $, $  \psi_1^+ (x)$,  and $\widetilde \psi_2 (x,t) $, $  \psi_2^+ (x)$.  Then by   Lemma~\ref{L1.2FLRW}   we obtain  
\begin{eqnarray*} 
&  &
 t^{ 3\ell/2-|\Im (m)|}\|\widetilde \psi_1 (t) -\widetilde \psi_2 (t)  \|_k\\
&   \leq  & 
c\|\widetilde \psi_1 (1) -\widetilde \psi_2 (1)\|_k  \\
&  \leq  & 
c  \|\psi_1 (1) - \psi_2 (1)\|_k+ \lim_{t \to  \infty}   \int_1^{ t} \left\|S(1,\tau )\left[  F(\psi_1 (x,\tau ))- F(\psi_2 (x,\tau )) \right] \,d \tau \right\|_k \\
& \leq  & 
c  \|\psi_1 (1) - \psi_2 (1)\|_k+ \lim_{t \to  \infty}   \int_1^{ t} c \tau ^{3\ell/2 +|\Im (m)|} \left\| F(\psi_1 (x,\tau ))- F(\psi_2 (x,\tau )) \right\|_k \,d \tau \\
&  \leq  & 
c  \|\psi_1 (1) - \psi_2 (1)\|_k\\
&  &
+ \lim_{t \to  \infty}   \int_1^{ t} c \tau ^{3\ell/2 +|\Im (m)|}\|  \psi_1 (x,\tau ))-  (\psi_2 (x,\tau )) \|_k \left(\|  \psi_1 (x,\tau ))\|_k^\alpha  + \|  \psi_2 (x,\tau ))\|_k^\alpha \right) \,d \tau  \\
&  \leq  & 
c  \|\psi_1 (1) - \psi_2 (1)\|_k\\
&  &
+c   \left( \sup_{\tau  \in [1,\infty) } \tau ^{3\ell/2-|\Im (m)|}
\|  \psi _1 (\tau)- \psi _2(\tau) \|_k \right)  \int_1^t 
\left(\|  \psi _1 (\tau)  \|_k^\alpha   
+\|    \psi _2 (\tau)\|_k^\alpha  \right)\,d \tau \\ 
& \leq  & 
c  \|\psi_1 (1) - \psi_2 (1)\|_k\\
&  &
+c   \left( \sup_{\tau  \in [1,\infty) } \tau ^{3\ell/2-|\Im (m)|}
\|  \psi _1 (\tau)- \psi _2(\tau) \|_k \right)   \left( \sup_{\tau  \in [1,\infty) } \tau ^{3\ell/2-|\Im (m)|}
\left(\|  \psi _1 (\tau)  \|_k   
+\|    \psi _2 (\tau)\|_k  \right)\right)^\alpha \,.
\end{eqnarray*} 
This completes the proof of the last statement of theorem. } 
\qed

\begin{corollary}
(i) For the matter dominated universe $\ell=2/3$,  the  condition of the theorem implies
\[
2|\Im (m)|+\alpha \left(|\Im (m)| -1 \right)  <-1 \Longleftrightarrow |\Im (m)|< 1 \quad and \quad \alpha >
1+ 3\frac{|\Im (m)|}{1-|\Im (m)| }\,.
\]
(ii) For the radiation dominated universe $\ell=1/2$, the  condition of the theorem implies
\[
2|\Im (m)|+\alpha \left(|\Im (m)| -\frac{3}{4}  \right)  <-1\Longleftrightarrow |\Im (m)|<\frac{3}{4} \quad and \quad \alpha >
1+ \frac{\frac{1}{4}+3|\Im (m)|}{\frac{3}{4}-|\Im (m)| }\,.
\]
\end{corollary}

\subsection{Finite speed propagation property}
\label{SS5.4}

In this section we prove that the dependence domain for the classical solution $u(t,x)$ at the point $(T,x_0)$  of the semilinear equation coincide with the one of the  corresponding linear equation. We present the proof here in order to make the paper self-contained. For { $(T,x_0) \in (1,\infty)\times \R^3$}  let  
\[
\Sigma _- {(T,x_0)}:= \left\{(t,x) \in { [1,T]}\times \R^3\,\Big|\, |x-x_0|=\phi (T)-\phi (t) \right\}
\]
be  a part of the backward   ``curved light cone'' (nullcone),  where  $\phi (t):= \frac{1}{1-\ell}t^{1-\ell}$ if $\ell\not=1$ and $\phi (t):= \ln(t)$ if $\ell =1$. Let  also     
\[
D_- {(T,x_0)}= \left\{(t,x) \in { [1,T]}\times \R^3\,\Big|\, |x-x_0|\leq \phi (T)-\phi (t)  \right\}
\]
be  the region defined in (\ref{Dpm}), whose       boundary contains   $\Sigma _- {(T,x_0)} $.
The outline of the proof of the next theorem is similar to one of  \cite{John}. 
\begin{theorem}
Let $\psi $ be a $C^1$ solution of the equation 
\begin{equation}
\label{5.4}
 \dsp 
\left(     \partial_0   +   \frac{1}{a(t)}\sum_{j=1,2,3}{\alpha }^j \partial_j+      \frac{3\dot a(t)}{2a(t)}{\mathbb I}_4     +i mt^{-1} {\gamma }^0 -iV (x,t)\right)\psi = { F(\psi )}  \,,
\end{equation}  
in the backward curved light cone $D_- (T,x_0)$ through $(T,x_0) \in(1,\infty)\times \R^3$. Assume that  the   potential $  V \in C([0,\infty) \times  \R^3 )  $ and  the nonlinear term $F(\psi    )\in C^1$ is   such that 
\begin{equation}
\label{5.F}
{ F(0 )=0}  \,.
\end{equation} 
If 
\begin{equation}
\label{5.IC}
\psi (x,1)=0 \quad \mbox{for all}\quad x \in D_- (T,x_0)\cap \{ t=1\}\,,
\end{equation}
then  $\psi $ vanishes in $ D_- (T,x_0)$.
\end{theorem}
\medskip

\ndt
{\bf Proof.} {  We consider the equation (\ref{5.4}) in the proper time 
\[
\tau := \phi (t)-\phi (1) \geq 0\ ,\quad \tau_0 := \phi (T)-\phi (1) >0,
\]
then $a(t)\partial_0= \partial_\tau $ and the equation reads 
\[
\left(     \partial_\tau    +   \sum_{j=1,2,3}{\alpha }^j \partial_j+  \ell b(\tau )    \frac{3 }{2 }{\mathbb I}_4     +i m  b(\tau )  {\gamma }^0 -ic(\tau ) V (x,t(\tau ))\right)\psi =c(\tau ){ F(\psi )}  \,,
\]
where
\[
\cases{ b(\tau ):=     \left[(1-\ell)\left(\tau +\phi (1)\right)\right]^{-1},\quad 
c(\tau ):=\left[(1-\ell)\left(\tau +\phi (1)\right)\right]^{\frac{\ell}{1-\ell}}\quad \mbox{\rm if} \quad \ell \not= 1,\cr
b(\tau ):=    1,\quad 
c(\tau ):=e^\tau\quad \mbox{\rm if} \quad\ell = 1.
}
\]
For $0 \leq s<\tau _{0}$ let
\[
\rho  (s, x)=\tau _{0}-\left[\left(\tau _{0}-s\right)^{2}+\tau _{0}^{-2}\left(2 \tau _{0} s-s^{2}\right)\left|x-x_{0}\right|^{2}\right]^{1 / 2} .
\]
Then
\[
\rho (0, x)=0, \quad \mbox{\rm and }  \quad  \lim _{s \rightarrow \tau _{0}} \rho (s, x)=\tau _{0}-\left|x-x_{0}\right|
\]
Define
\[
R_{s}=\left\{(\tau , x)\,|\, 0 \leq \tau  \leq \rho (s, x),\left|x-x_{0}\right|<\tau _{0}-\tau \right\},
\]
then in the proper time  the interior of $D_- {(T,x_0)} $ is
\[
\Lambda_{-} {\left(\tau _{0}, x_{0}\right)}:= \left\{(\tau , x) \in\left[0,\tau _{0}\right) \times {\mathbb R}^{n}\,|\,\left|x-x_{0}\right|<\tau _{0}-\tau \right\} =\bigcup_{0 \leq s<\tau _{0}} R_{s}\,.
\]
}
{ 
We fix $  s_{0} $  such that $0 < s_{0}<\tau _{0}$, then with some number $\theta\left(s_{0}\right) $ we have  
\[
\left|\nabla_{x} \rho (s, x)\right|=\frac{\tau _{0}^{-2}\left(2 \tau _{0} s-s^{2}\right)\left|x-x_{0}\right|}{\left[\left(\tau _{0}-s\right)^{2}+\tau _{0}^{-2}\left(2 \tau _{0} s-s^{2}\right)\left|x-x_{0}\right|^{2}\right]^{1 / 2}} \leq \theta\left(s_{0}\right)<1
\]
if $0 \leq s \leq s_{0}<\tau _{0}$. 
Next we consider the surface defined by
\[
\Lambda_{s}:=\left\{(\tau , x)\,|\, \tau =\rho (s, x),\left|x-x_{0}\right|<\tau _{0}\right\}\,.
\]
Note that the outward unit normal at $(\rho (s, x), x) \in \Lambda_{s}$ is
\[
\frac{ \left(1,-\nabla_{x} \rho \right) }{ \sqrt{1+\left|\nabla_{x} \rho \right|^{2}} }.
\]
To apply the energy method   to the equation (\ref{5.4}),
we use the identity 
\begin{eqnarray*}
&  &
  \dsp 
   \partial_\tau  |\psi|^2    +\sum_{j=1,2,3} \partial_j \left(  \psi ^*  \alpha ^j  \psi \right) +3 \ell b(\tau )|\psi|^2    -2\Im (m) b(\tau )\psi ^*  {\gamma }^0 \psi+2\psi ^*   c(\tau )\Im (V(x,t(\tau)))  \psi \\
& = &
 c(\tau )(\psi ^* F(\psi ) +   F^*(\psi ) \psi )\,.  
\end{eqnarray*} 
The divergence theorem and the vanishing initial data (\ref{5.IC}) imply
\begin{eqnarray*}
 &  &
\int_{R_s}  \dsp 
  \left(  \partial_\tau  |\psi|^2    +\sum_{j=1,2,3} \partial_j \left(  \psi ^*  \alpha ^j \psi \right) \right) d\tau  \,dx \\
& =  &
\int_{\Lambda _s}  \dsp 
  \left(  |\psi|^2    -\sum_{j=1,2,3}  \left( \psi ^*  \alpha ^j  \psi \right)\partial_j \rho    (s,x)  \right) \frac{1}{\sqrt{1+|\nabla_x \rho  (s,x)|^2}}\,d \sigma\,. 
 \end{eqnarray*} 
} 
{ Therefore  
\begin{eqnarray}
\label{5.5}
& &
\int_{\Lambda _s}  \dsp 
  \left(  |\psi|^2    -\sum_{j=1,2,3}  \left( \psi ^*  \alpha ^j  \psi \right)\partial_j \rho  (s,x)  \right) \frac{1}{\sqrt{1+|\nabla_x \rho   (s,x)|^2}}\,d \sigma\nonumber \\
&  &
+    \int_{R_s}  \dsp 
  \left( 3b(\tau)|\psi|^2    -2\Im (m)b(\tau)\psi ^*  {\gamma }^0 \psi +2c(\tau)\psi ^*  \Im (V(x,t(\tau)))  \psi\right) d\tau  \,dx  \nonumber \\
& = &
\int_{R_s} c(\tau) \left(  \psi ^* F(\psi  ) + F^*(\psi  ) \psi  \right) d\tau  \,dx\,.  
\end{eqnarray}
The hermitian matrix
\[
A={\mathbb I}_4- \sum_{j=1,2,3}    \alpha ^j  a_j,\quad a_j:= \partial_j \rho   (s,x) ,\quad j=1,2,3,
\]
 for all $s \in [0,s_0]$ has two double positive eigenvalues  
$1-\sqrt{a_1^2+a_2^2+a_3^2}>0$  and $1+\sqrt{a_1^2+a_2^2+a_3^2}$, since $\left|\nabla_{x} \rho (s, x)\right| \leq \theta\left(s_{0}\right)<1$. 
Hence, there is $\delta (s_0)>0 $ such that
\[
|\psi|^2    -\sum_{j=1,2,3}  \left( \psi ^*  \alpha ^j  \psi \right)\partial_j \rho  (s,x) \geq  
\delta (s_0) | \psi |^2 \quad  \mbox{\rm for all}\quad s \in [0,s_0]\,.
\]
The equation (\ref{5.5}) and condition (\ref{5.F}) yield
\begin{eqnarray*}
&  &
\int_{\Lambda _s}  \dsp 
   |\psi|^2   \frac{1}{\sqrt{1+|\nabla_x \rho   (s,x)|^2}}\,d \sigma \\
& \leq  &
C(s_0) \int_0^s \int_{\Lambda _\lambda }  \dsp 
 |\psi|^2    \frac{\partial_\lambda\rho    (\lambda ,x)}{\sqrt{1+|\nabla_x \rho   (\lambda ,x)|^2}}\,d \sigma\,  d\lambda  
\quad  \mbox{\rm for all}\quad s  \in [0,s_0] \,.
\end{eqnarray*}
}
{ 
Now we set
\[
I(s)=\int_{\Lambda _s}  \dsp 
|\psi|^2   \frac{1}{\sqrt{1+|\nabla_x \rho  |^2}}\,d \sigma \,,
\]
and derive 
\[
I(s)
\leq 
C(s_0) \left(\max_{ 0\leq s\leq s_0 }\left|\partial_s \rho    (s ,x) \right|\right) \int_0^s  I(\lambda ) \, d\lambda  
\quad  \mbox{\rm for all}\quad  0\leq s \leq s_0\,.
\]   
It remains to apply Gronwall's inequality. 
\qed
}

\section{Lochak-Majorana condition in FLRW spacetime}
\label{S2}
\setcounter{equation}{0}

Consider the Dirac equation  with the matrix-valued potential function 
\begin{equation}
\label{A}
A(x,t) = \alpha (x,t){\mathbb I}_4 +i\beta (x,t) \gamma ^5, 
\end{equation}
 where { $\alpha (x,t),  \beta (x,t)$ are any real-valued continuous} functions.

\begin{lemma}
\label{L2.1}
For the solution $\psi \in C^1([1,\infty); L^2({\mathbb R}^3))\cap C^0([1,\infty); H_{(1)}({\mathbb R}^3))$ of the  Dirac equation
\begin{equation}
\label{DE_LM}
  \dsp 
\left(   i {\gamma }^0    \partial_0   +i t^{-\ell}  \sum_{j=1,2,3} {\gamma }^j  \partial_{j}+i \frac{3\ell  }{2} t^{-1}    {\gamma }^0     -m t^{-1}{\mathbb I}_4  + {\gamma }^0 V(x,t)\right)\psi=-A \psi  \,,
\end{equation}
the following energy identity holds
\begin{eqnarray}
  \dsp 
\|\psi (\cdot ,t)\|_{L^2({\mathbb R}^3)}^2 
 & = &
t^{-3\ell} 
\|\psi  (\cdot ,1)\|_{L^2({\mathbb R}^3)}^2  +
2\Im (m)t^{-3\ell}\int_1^t s^{3\ell-1}  
\int_{{\mathbb R}^3}\psi ^* (x ,s) {\gamma }^0 \psi (x ,s)\, dx\,ds \nonumber \\
\label{2.3}
&  &
- 2 t^{-3\ell}\int_1^t s^{ 3\ell} 
\int_{{\mathbb R}^3} \psi ^* (x ,s) \Im (V(x,s))  \psi(x ,s)  \, dx \,ds
\quad \mbox{for all}\quad t\geq 1.
\end{eqnarray}
\end{lemma}
\medskip

\ndt
{\bf Proof.} 
We write the equation (\ref{DE_LM}) in the equivalent form
\begin{eqnarray*}
  \dsp 
\left(     \partial_0   +   t^{-\ell}\sum_{j=1,2,3}{\alpha }^j \partial_j+      \frac{3\ell  }{2} t^{-1} {\mathbb I}_4     +i mt^{-1} {\gamma }^0 -iV (x.t)\right)\psi=  i{\gamma }^0A \psi  \,.
\end{eqnarray*}
Hence
\[
  \dsp 
   \partial_0 |\psi|^2    +\sum_{j=1,2,3} \partial_j \left(  t^{-\ell}\psi ^*  \alpha ^j  \psi \right) +3\ell t^{-1}|\psi|^2    -2\Im (m)t^{-1}\psi ^*  {\gamma }^0 \psi+2\psi ^*  \Im (V(x,t))  \psi \\
  =  
\psi ^* i{\gamma }^0A \psi + (i{\gamma }^0A \psi)^*  \psi \,.  
\]  
The terms with  the potential  $A$ can be easily evaluated: 
\begin{eqnarray*}
 \psi ^* i{\gamma }^0A \psi + (i{\gamma }^0A \psi)^*  \psi  
&  = &0\,.
\end{eqnarray*}
Therefore, 
\begin{eqnarray*}
&  &
  \dsp 
   \partial_0 |\psi|^2   =- \sum_{j=1,2,3} \partial_j \left( t^{-\ell}\psi ^*  \alpha ^j  \psi \right) -3\ell t^{-1}|\psi|^2   + 2\Im (m)t^{-1}\psi ^*  {\gamma }^0 \psi - 2\psi ^*  \Im (V(x,t))  \psi  \,.  
\end{eqnarray*}
The finite propagation speed property implies
\begin{eqnarray*}
&  &
  \dsp 
\frac{d}{dt} \|\psi \|_{L^2({\mathbb R}^3)}^2  =-3\ell t^{-1}\|\psi \|_{L^2({\mathbb R}^3)}^2   + 2\Im (m)t^{-1}\int_{{\mathbb R}^3}\psi ^*  {\gamma }^0 \psi \, dx- 2\int_{{\mathbb R}^3} \psi ^*  \Im (V(x,t))  \psi  \, dx 
\end{eqnarray*}
and, consequently, 
\[
  \dsp 
\frac{d}{dt} \left( t^{3\ell} \|\psi \|_{L^2({\mathbb R}^3)}^2  \right) 
  =  
2\Im (m)t^{-1} t^{3\ell}\int_{{\mathbb R}^3}\psi ^*  {\gamma }^0 \psi \, dx- 2 t^{3\ell}\int_{{\mathbb R}^3} \psi ^*  \Im (V(x,t))  \psi  \, dx \,.
\]
It follows (\ref{2.3}). 
Lemma is proved. \qed

\begin{lemma}
\label{L3.2}
Assume that
\begin{equation}
\label{Vgamma2}
 V^T(x,t)\gamma ^2  +  \gamma ^2V(x,t)=0\,.
\end{equation}
For the solution $\psi \in C^1([0,\infty); L^2({\mathbb R}^3))\cap C^0([0,\infty); H_{(1)}({\mathbb R}^3))$ of the  Dirac 
equation (\ref{DE_LM}) we have
\[
\int_{{\mathbb R}^3}  \psi^T (x,t) \gamma ^2\psi (x,t) \,dx
 =
t^{-3\ell} 
 \int_{{\mathbb R}^3}  \psi^T (x,1) \gamma ^2\psi (x,1) \,dx  \,.
\]
In particular, for $z \in {\mathbb C}$, $|z|=1$,   
\begin{eqnarray*}
\int_{{\mathbb R}^3} 2\Re \Big( \ol{z}\psi^T (x,t) \gamma ^2\psi (x,t)\Big) \,dx 
 & = &
t^{-3\ell}  \int_{{\mathbb R}^3} 2\Re \Big(  \ol{z}\psi^T (x,1) \gamma ^2\psi (x,1) \Big) \,dx\,. 
\end{eqnarray*}
\end{lemma}
\medskip

\noindent
{\bf Proof.} We multiply the equation (\ref{DE_LM}) by { $-i\psi^T\gamma ^2 \gamma ^0$}
and obtain 
\[
 \dsp 
\psi^T\gamma ^2  \partial_0  \psi   
  = 
- t^{-\ell} \sum_{j=1,2,3}\psi^T\gamma ^2{\gamma }^0{\gamma }^j  \partial_{x_j} \psi  
 -\frac{3\ell  }{2 }t^{-1} \psi^T\gamma ^2     \psi   -im t^{-1}\psi^T\gamma ^2{\gamma }^0 \psi+\psi^T\gamma ^2iV(x,t)\psi +i\psi^T\gamma ^2{ {\gamma }^0A} \psi  .
\]
Here
\[
\gamma ^2\gamma ^0 
= \left(
\begin{array}{cccc}
 {\mathbb O}_2 & -\sigma _2  \\
 -\sigma _2 & {\mathbb O}_2  \\
\end{array}
\right), \quad \gamma ^2\gamma ^0\gamma ^1=
\left(
\begin{array}{cccc}
-i\sigma _3   & {\mathbb O}_2 \\
 {\mathbb O}_2 & i\sigma _3  \\
\end{array}
\right)
, \quad\gamma ^2\gamma ^0\gamma ^2=\gamma ^0,\quad \gamma ^2\gamma ^0\gamma ^3
=\left(
\begin{array}{cccc}
i\sigma _1   & {\mathbb O}_2 \\
 {\mathbb O}_2 & -i\sigma _1  \\
\end{array}
\right),
\]
\[
\gamma ^0\gamma ^5=
\left(
\begin{array}{cccc}
 {\mathbb O}_2 & -{\mathbb I}_2  \\
 {\mathbb I}_2 & {\mathbb O}_2  \\
\end{array}
\right)=-\gamma ^5\gamma ^0,   \quad \gamma ^5\gamma ^2+\gamma ^2\gamma ^5=0\,.
\]
Now 
\begin{eqnarray*}
  \dsp 
 \partial_0  \psi^T  =- t^{-\ell}   \sum_{j=1,2,3} \partial_j\psi ^T ({\gamma ^j})^T \gamma ^0 -    \frac{3\ell  }{2 }t^{-1}\psi ^T    -i mt^{-1} \psi^T {\gamma }^0 +  i\psi ^TV^T(x,t)+ i\psi^T A {\gamma }^0  \,,
\end{eqnarray*}
and, consequently,
\begin{eqnarray*}
 \partial_0  (\psi^T\gamma ^2\psi)
 & = &
\Bigg(-  t^{-\ell}\sum_{j=1,2,3} \partial_j\psi ^T ({\gamma ^j})^T \gamma ^0 -    \frac{3\ell  }{2 }t^{-1}\psi ^T    -i mt^{-1} \psi ^T{\gamma }^0 + i\psi^T A {\gamma }^0  \Bigg)\gamma ^2 \psi\\
&  &
- t^{-\ell} \sum_{j=1,2,3}\psi^T\gamma ^2{\gamma }^0{\gamma }^j  \partial_{j} \psi  -\frac{3\ell  }{2 }t^{-1} \psi^T\gamma ^2     \psi   -im t^{-1}\psi^T\gamma ^2{\gamma }^0 \psi\\
&  &
+ i\psi ^T \left( V^T(x,t)\gamma ^2  +  \gamma ^2V(x,t)\right) \psi+i\psi^T\gamma ^2{ {\gamma }^0A} \psi \,.
\end{eqnarray*}
Consider now the terms with the potential $A$:
\begin{eqnarray}
\Big( i\psi^T  A{\gamma }^0  \Big)\gamma ^2 \psi
+ \psi^T\gamma ^2   \Big( i{ {\gamma }^0A} \psi \Big)
 & = &
i\psi^T  \Bigg(\alpha I +i\beta \gamma ^5\Bigg){\gamma }^0 \gamma ^2 \psi
+   i\psi^T\gamma ^2  { {\gamma }^0\Bigg(\alpha I +i\beta \gamma ^5 \Bigg)} \psi \nonumber \\
 & = &
i\psi^T  \Bigg(\alpha  {\gamma }^0\gamma ^2+i\beta \gamma ^5{\gamma }^0 \gamma ^2
+   {   \alpha \gamma ^2 {\gamma }^0 +i\beta  \gamma ^2  {\gamma }^0\gamma ^5 }\Bigg)\psi \nonumber \\
\label{2.4}
 & = &
-\psi^T \beta  \Big(\gamma ^5{\gamma }^0 \gamma ^2
+\gamma ^2 {  {\gamma }^0\gamma ^5 }\Big)\psi=0\,.
\end{eqnarray}
Thus, from (\ref{Vgamma2}) and (\ref{2.4}) we obtain
\begin{eqnarray*} 
\partial_0  (\psi^T\gamma ^2\psi)
 & =  &
  -3\ell t^{-1} \psi^T\gamma ^2     \psi -  t^{-\ell}\sum_{j=1,2,3} \partial_j\psi ^T ({\gamma ^j})^T \gamma ^0 \gamma ^2 \psi
- t^{-\ell}\sum_{j=1,2,3}\psi^T\gamma ^2{\gamma }^0{\gamma }^j  \partial_{j} \psi \\
 & =  &
  -3\ell t^{-1}  \psi^T\gamma ^2     \psi -  t^{-\ell}\sum_{j=1,2,3} \Big\{ (\partial_j\psi ^T )({\gamma ^j})^T \gamma ^0 \gamma ^2 \psi
+ \psi^T\gamma ^2{\gamma }^0{\gamma }^j  \partial_{ j} \psi \Big\} \,.
\end{eqnarray*}
For the sum in the last equation we have
\begin{eqnarray*}
({\gamma }^1)^T{\gamma }^0 \gamma ^2
= \gamma ^2{\gamma }^0{\gamma }^1\,,\quad
({\gamma }^2)^T{\gamma }^0 \gamma ^2 
= \gamma ^2{\gamma }^0  \gamma ^2\,,\quad
({\gamma }^3)^T{\gamma }^0 \gamma ^2
= \gamma ^2{\gamma }^0{\gamma }^3\,.
\end{eqnarray*}
It follows 
\begin{eqnarray*} 
\partial_t  (\psi^T\gamma ^2\psi)
 & =  &
  -3\ell t^{-1}\psi^T\gamma ^2     \psi -  t^{-\ell}\sum_{j=1,2,3} \partial_j \Big(\psi ^T    \gamma ^2 \gamma ^0{\gamma ^j}\psi\Big) 
\end{eqnarray*}
and, consequently, 
\[
\frac{d}{dt} \int_{{\mathbb R}^3}  \psi^T (x,t) \gamma ^2\psi (x,t) \,dx
 =
 - 3\ell t^{-1} \int_{{\mathbb R}^3} \psi^T  (x,t) \gamma ^2 \psi (x,t) \,dx\,.
\]
Thus, the first statement of the lemma is proved. To prove the second statement, we use 
$\ol{\gamma _2}=  -\gamma _2$ and  recall the formula from  \cite{Bachelot}:
\begin{eqnarray}
\label{3.17}
&  &
|\psi -z\gamma ^2 \ol{\psi }|^2= 2|\psi |^2 +2 \Re (\ol{z} \psi^T \gamma ^2 \psi) \,,
\end{eqnarray}
where $z \in {\mathbb C}$,   $|z|=1$. 
The lemma is proved. \qed

\begin{lemma}
Assume that    { $V^*(x,t)=V(x,t)$ and} 
(\ref{Vgamma2}) if fulfilled.
For the solution $\psi \in C^1([1,\infty); L^2({\mathbb R}^3))\cap C^0([1,\infty); H_{(1)}({\mathbb R}^3))$ of the  Dirac 
equation (\ref{DE_LM}) one has
\begin{eqnarray*}
\int_{{\mathbb R}^3} |\psi (x,t) -z\gamma ^2 \ol{\psi  (x,t)}|^2\,dx 
& =  &
t^{-3\ell} \left( \int_{{\mathbb R}^3} |\psi (x,1) -z\gamma ^2 \ol{\psi  (x,1)}|^2\,dx
\right) \\
&  &
+  4\Im (m)t^{-3\ell}\int_1^t  s^{3\ell-1}\int_{{\mathbb R}^3}\psi ^* (x ,s) {\gamma }^0 \psi (x ,s)\, dx\,ds \,.
\end{eqnarray*}
\end{lemma}
\medskip

\noindent
{\bf Proof.} From Lemma~\ref{L2.1}, Lemma~\ref{L3.2}, and (\ref{3.17})
\begin{eqnarray*}
\int_{{\mathbb R}^3} |\psi (x,t) -z\gamma ^2 \ol{\psi  (x,t)}|^2\,dx 
& =  &
\int_{{\mathbb R}^3} \Big(2|\psi |^2 +2 \Re (\ol{z} \psi^T \gamma ^2 \psi)\Big)\,dx \\
& =  &
2 \left( t^{-3\ell}  \|\psi  (\cdot ,1)\|_{L^2({\mathbb R}^3)}^2  +
2\Im (m)t^{-3\ell} \int_1^t s^{3\ell-1} \int_{{\mathbb R}^3}\psi ^* (x ,s) {\gamma }^0 \psi (x ,s)\, dx\,ds\right) \\
&  &
+ t^{-3\ell}  \int_{{\mathbb R}^3} 2\Re \Big(  \ol{z}\psi^T (x,1) \gamma ^2\psi (x,1) \Big) \,dx\,.
\end{eqnarray*}
That is,
\begin{eqnarray*}
\int_{{\mathbb R}^3} |\psi (x,t) -z\gamma ^2 \ol{\psi  (x,t)}|^2\,dx 
& =  &
2  t^{-3\ell} \left( \|\psi  (\cdot ,1)\|_{L^2({\mathbb R}^3)}^2  + \int_{{\mathbb R}^3} \Re \Big(  \ol{z}\psi^T (x,1) \gamma ^2\psi (x,1) \Big) \,dx
\right) \\
&  &
+  4\Im (m)t^{-3\ell} \int_1^t s^{3\ell-1}\int_{{\mathbb R}^3}\psi ^* (x ,s) {\gamma }^0 \psi (x ,s)\, dx\,ds \,.
\end{eqnarray*}
The lemma is proved. \qed

The last statement of the next corollary contains   (\ref{LMC}) written as follows:
\[
\rho ^2 (\psi )
 = 
(|\psi _1|^2+|\psi _2|^2-|\psi _3|^2-|\psi _4|^2)^2+ (2\Im (\psi _1\ol{\psi _3})+2\Im (\psi _2\ol{\psi _4}))^2\,.
\]
\begin{corollary}
\label{C2.5}
Assume that  
(\ref{Vgamma2}) if fulfilled.
(i) If $\Im (m)=0$, then
\begin{eqnarray*}
\int_{{\mathbb R}^3} |\psi (x,t) -z\gamma ^2 \ol{\psi  (x,t)}|^2\,dx 
& =  &
  t^{-3\ell}\left( \int_{{\mathbb R}^3} |\psi (x,1) -z\gamma ^2 \ol{\psi  (x,1)}|^2\,dx
\right)  \,.
\end{eqnarray*}
(ii) If $\psi (x,1) -z\gamma ^2 \ol{\psi  (x,1)}=0$ and $\Im (m) \not=0$, then
\begin{eqnarray*}
\int_{{\mathbb R}^3} |\psi (x,t) -z\gamma ^2 \ol{\psi  (x,t)}|^2\,dx 
& =  &
  4\Im (m) t^{-3\ell}\int_1^t  s^{3\ell-1}\int_{{\mathbb R}^3}\psi ^* (x ,s) {\gamma }^0 \psi (x ,s)\, dx\,ds \,.
\end{eqnarray*}
(iii) if $\psi (x,1) -z\gamma ^2 \ol{\psi  (x,1)}=0$ , then
\begin{eqnarray*}
\int_{{\mathbb R}^3} |\psi (x,t) -z\gamma ^2 \ol{\psi  (x,t)}|^2\,dx 
& \leq   & 
 4|\Im (m)| t^{-3\ell}\int_1^t s^{3\ell-1}\int_{{\mathbb R}^3}\rho  (x ,s)\, dx\,ds  \,. 
\end{eqnarray*}
\end{corollary}

\ndt
{\bf Proof.} It is evident. \qed

Proposition~I.1~\cite{Bachelot} implies  
that $\psi -z\gamma ^2 \ol{\psi }=0$  is equivalent to $\rho ^2 (\psi )=0 $ if  $z \in {\mathbb C}$, $|z|=1$.

\section{Large amplitude global solution. Proof of Theorem~\ref{T4.1}}
\label{S3}
\setcounter{equation}{0}

We are going to prove the existence of the global solution by the successive approximations. We follow approach developed by Bachelot~\cite{Bachelot} that appeals to the estimates obtained by the replacing the generators of the Poincar\'e  group with the Fermi operators. Energy estimates obtained in subsection~\ref{SSEE}
allows  us  to expand  the result from \cite{Bachelot} to to the Dirac equation in the FLRW spacetime. 

\medskip

  Let $\Psi $ be the solution of  
\begin{equation}
\label{B31}  
\cases{ \dsp\left( i {\gamma }^0    \partial_0   +i t^{-\ell}\sum_{j=1,2,3}{\gamma }^j  \partial_{j} +i \frac{3\ell  }{2 }  t^{-1}  {\gamma }^0     -m t^{-1}{\mathbb I}_4  +i \gamma ^0V(x,t) \right)\Psi (x,t)  =0  \,,\cr
\Psi (x,1)=\Psi _0 (x)\,.
}
\end{equation}
Here $ V^*(x,t)=V(x,t)$. The finite propagation speed property  implies that  for every $t>1$ the support of $\Psi =\Psi (x,t) \in C^\infty({\mathbb R}^3\times (1,\infty)  ; {\mathbb C}^4)$ is in  some compact subset 
of ${\mathbb R}^3$. Then according to Lemma~\ref{L1.2FLRW}
\begin{equation}
\label{4.26FLRW}
\|\Psi (t) \|_k
\leq 
c t^{-3\ell/2+|\Im (m)|}\|\Psi _0 \|_k  \quad \mbox{\rm for all}\quad t\geq 1\,. 
\end{equation}
 We look for the solution of (\ref{NDE_CP}) in the form
\begin{equation}
\label{B33}
\psi =\Psi +\chi \,,
\end{equation}
where $\Psi  $ solves (\ref{B31}). 

We consider  the nonlinear term.
The Lochak-Majorana condition (\ref{LMC}) 
 is equivalent to $\rho (\Psi_0 )=0$ (see \cite{Bachelot}).  First  we note, that according to Corollary~\ref{C2.5},   (\ref{Falphabeta}),    and (\ref{LMC}),  
if $\Im (m)=0$, then 
\begin{eqnarray*}
\dsp F\left( \Psi ^*  (x,t) \gamma ^0 \Psi   (x,t) \,,\,\Psi ^*  (x,t)\gamma ^0\gamma ^5 \Psi    (x,t) \right)
& = &
0 \quad \mbox{\rm for all}\quad t \geq 1,\quad x \in {\mathbb R}^3 \, .
\end{eqnarray*}
Further, (\ref{F}), 
where $\alpha $ and $\beta $ are real-valued functions,  
can be written as follows
\[
F\left(\xi ,\eta  \right)= \alpha \left(\xi ,\eta  \right)I +i\beta \left(\xi ,\eta  \right) \gamma ^5,\quad \xi = \ol{\psi  }^T\gamma ^0 \psi \in {\mathbb R},\quad
\eta =\ol{\psi  }^T\gamma ^0\gamma ^5 \psi \in {\mathbb R}\,.
\]
It is evident that with the solution $\Psi   (x,t)  $ we can write
\begin{eqnarray*}
|\ol{(\Psi (x,t)+\chi (x,t))}^T\gamma ^0 (\Psi(x,t) +\chi (x,t))-\ol{ \Psi (x,t) }^T\gamma ^0  \Psi(x,t)| 
& \!\! \leq \!\!  &
C{ |\chi (x,t)| \big(|\chi (x,t)|+|\Psi  (x,t)|\big)},\\
|\ol{(\Psi (x,t)+\chi (x,t))  }^T\gamma ^0\gamma ^5 (\Psi (x,t)+\chi (x,t))- \ol{(\Psi  )  }^T\gamma ^0\gamma ^5 (\Psi  (x,t))|
& \!\! \leq  \!\! &
C { |\chi (x,t)|\big(|\chi (x,t)|+|\Psi  (x,t)|\big)} 
\end{eqnarray*}
and, consequently, 
\begin{eqnarray*}
 \alpha \left( \psi^*\gamma ^0 \psi\,,\, \psi^*\gamma ^0\gamma ^5 \psi \right)\psi 
& = &
  \alpha \left( (\Psi +\chi )^*\gamma ^0 (\Psi +\chi )\,,\, (\Psi +\chi )^*\gamma ^0\gamma ^5 (\Psi +\chi ) \right)\psi \\
& = &
  \alpha \left(  \Psi^*\gamma ^0 \Psi  \,,\, \Psi^*\gamma ^0\gamma ^5 \Psi  \right)\psi +\alpha_1 \left(\chi ,\Psi \right)  \\
& = &
\alpha_1 \left(\chi ,\Psi \right)\,,\\
 \beta  \left(\psi^*\gamma ^0 \psi\,,\,\psi^*\gamma ^0\gamma ^5 \psi \right)\psi 
& = &
  \beta \left( (\Psi +\chi )^*\gamma ^0 (\Psi +\chi )\,,\, (\Psi +\chi )^*\gamma ^0\gamma ^5 (\Psi +\chi ) \right)\psi \\
& = &
  \beta  \left( \Psi^*\gamma ^0 \Psi  \,,\, \Psi^*\gamma ^0\gamma ^5 \Psi  \right)\psi +\beta _1 \left(\chi ,\Psi \right)  \\
& = &
\beta _1 \left(\chi ,\Psi \right)\,,
\end{eqnarray*}
where  $ \alpha_1 ,\beta _1  \in C^\infty({\mathbb C}^8; {\mathbb C}^4)$, and, as the functions of $\chi  $,
\begin{eqnarray}
\label{B35a}
|\alpha _1\left(\chi \,,\,\Psi  (x,t)\right)| 
& = &
O\Big( { |\chi |\big(|\chi |+|\Psi (x,t) |\big) } \Big)\quad \mbox{\rm as} \quad |\chi| \to 0\,, \\
\label{B35b}
|\beta  _1\left(\chi \,,\,\Psi  (x,t)\right)| 
& = &
O\Big( { |\chi |\big(|\chi |+|\Psi (x,t) |\big) } \Big)\quad \mbox{\rm as} \quad |\chi| \to 0\,.
\end{eqnarray}

Further,
the Cauchy problem becomes 
\begin{eqnarray}
\label{B34}
\cases{ \dsp\left(i {\gamma }^0    \partial_0   +i t^{-\ell}\sum_{j=1,2,3}{\gamma }^j  \partial_{j} +i \frac{3\ell  }{2 }  t^{-1}   {\gamma }^0     -m t^{-1}{\mathbb I}_4  +i \gamma ^0V(x,t) \right)\chi   =f_1\left(\chi \,,\,\Psi  \right)   \,,\cr
\chi  (x,1)=\varepsilon \chi _0 (x)\,,
}
\end{eqnarray}
where, in view of (\ref{B31}), we have denoted
\begin{eqnarray*}
f_1\left(\chi \,,\,\Psi  \right)
& := &
  F\left((\Psi  +\chi )^*\gamma ^0 (\Psi  +\chi  )\,,\,(\Psi  +\chi )^*\gamma ^0\gamma ^5 (\Psi  +\chi  ) \right)  \Psi  \\
&   &
 + F\left((\Psi  +\chi )^*\gamma ^0 (\Psi  +\chi  )\,,\,(\Psi  +\chi )^*\gamma ^0\gamma ^5 (\Psi  +\chi  ) \right)  \chi    
\end{eqnarray*}
that can be rewritten as $ f_1\left(\chi \,,\,\Psi  \right)=-A \chi $ similar to (\ref{A}). 
Hence,  $ f_1  \in C^\infty({\mathbb C}^8; {\mathbb C}^4)$, while (\ref{B35a}) and (\ref{B35b}) imply
\[
|f_1\left(\chi \,,\,\Psi  (x,t)\right)| =O\Big( |\chi |\big(|\chi |+|\Psi (x,t) |\big)^2\Big)\quad \mbox{\rm as} \quad |\chi| \to 0\,.
\]
We look for the function $\chi  $ as a limit of the sequence  $\{\chi ^{(k)}\}_1^\infty $ that is defined  by
\[ 
\cases{ \dsp\left(  i {\gamma }^0    \partial_0   +it^{-\ell}\sum_{j=1,2,3}{\gamma }^j  \partial_{j} +i \frac{3\ell  }{2 }  t^{-1}   {\gamma }^0      -mt^{-1}{\mathbb I}_4 +i \gamma ^0V(t) \right)\chi^{(k)}   =f_1\left(\chi^{(k-1)} \,,\,\Psi  \right)   \,,\cr
\chi ^{(k)} (x,1)=\varepsilon \chi _0 (x), \quad k=1,2,\ldots\,,
}
\]
and $\chi ^{(0)}(x,t)\equiv 0 $. The finite propagation speed property implies that for every $t>0$ the supports of the functions $\chi^{(k)} =\chi^{(k)} (x,t) $ are in the same compact subset 
of ${\mathbb R}^3$ for all $k \geq 0$. Lemmas~\ref{L1.2FLRW}, \ref{L2.1}  and 
\begin{eqnarray*}
\| \chi ^{(1)}(t)\|_s 
& \leq  & 
c t^{-3\ell/2} \|\chi ^{(1)}(1) \|_s +
 c    t^{-3\ell/2} \int_1^t \tau^{3\ell/2}  \|f_1\left(\chi^{(0)}(\tau ) \,,\,\Psi (\tau ) \right)\|_s \,d\tau    \,,\quad  t\geq 1\,,
\end{eqnarray*}
imply 
\begin{eqnarray*}
\| \chi ^{(1)}(t)\|_s 
& \leq  & 
c t^{-3\ell/2} \varepsilon \|\chi _0  \|_s     \,,\quad  t\geq 1\,.
\end{eqnarray*} 
Corollary~6.4.5~\cite{Hor} and Lemma~\ref{L1.2FLRW} imply for $k=2,3,\ldots $ 
the following estimate
\begin{eqnarray}
\label{3.8}
\| \chi ^{(k)}(t)\|_s 
& \leq  & 
c   t^{-3\ell/2}\varepsilon \|\chi _0 (x) \|_s +
 c    t^{-3\ell/2} \int_1^t \tau^{3\ell/2}  \| f_1\left(\chi^{(k-1)} \,,\,\Psi  \right) (\tau )\|_s \,d\tau   \\ 
& \leq  & 
c  t^{-3\ell/2} \varepsilon \|\chi _0 (x) \|_s \nonumber 
 + c   t^{-3\ell/2}  \int_1^t \tau^{3\ell/2}   \| \chi^{(k-1)} (\tau )\|_s \Big( | \chi^{(k-1)} (\tau )|_{\left[\frac{s}{2}\right]}
+ |\Psi  (\tau )|_{\left[\frac{s}{2}\right]} \Big)^2   \,d\tau  ,\,\,  t\geq 1, \nonumber 
\end{eqnarray}
where  $\left[\frac{s}{2}\right] $ is the integer part of $\frac{s}{2} $, while
\[
|\Psi  (t )|_{s}:=\sup_{|\alpha| \leq s}  \|\partial_x^\alpha \Psi  (x, t )\|_{L^\infty({\mathbb R}^3)}.
\] 
Now we apply Sobolev embedding theorem and Lemma~\ref{L1.2FLRW}  to the function $ \Psi $:  
\begin{eqnarray*}
 | \Psi  (t)  |_{\left[\frac{s}{2}\right]} \leq \| \Psi  (t) \|_{s} 
& \leq  &
c  t^{-3\ell/2} \| \Psi_0  \|_{s}\quad \mbox{\rm for all}\quad t\geq 1\,.
\end{eqnarray*}
For the given $s$ and $n$ we define 
\[
a_n(t):=\sup_{{1\leq \tau \leq t,\,\,0\leq k\leq n}}\tau ^{3\ell/2}\| \chi^{(k)}  (\tau ) \|_{s}\,,\qquad A_n:=\sup_{t \in [1,\infty)} a_n(t)\,.
\]
Then  for { $s \geq 3$} by (\ref{3.8}) we derive
\begin{eqnarray*}
\| \chi ^{(k)}(t)\|_s 
& \leq  & 
c  t^{-3\ell/2} \varepsilon \|\chi _0 (x) \|_s \\
&  &
 + c   t^{-3\ell/2}  \int_1^t \tau ^{3\ell/2}  \| \chi^{(k-1)} (\tau )\|_s \Big( \| \chi^{(k-1)} (\tau )\|_s
+ \|\Psi  (\tau )\|_s \Big)^2   \,d\tau   \,,  \\
& \leq  & 
c    t^{-3\ell/2}   \Bigg\{ \varepsilon \|\chi _0 (x) \|_s \\
&  &
 +  \int_1^t \tau ^{-3\ell/2} \| \chi^{(k-1)} (\tau )\|_s \Big( \tau ^{3\ell/2}\| \chi^{(k-1)} (\tau )\|_s
+ \tau ^{3\ell/2}\|\Psi  (\tau )\|_s \Big)^2   \,d\tau  \Bigg\},\quad  t\geq 1\,,
\end{eqnarray*}
and 
\begin{eqnarray*} 
a_n(t) 
& \leq  & 
c   \varepsilon \|\chi _0 (x) \|_s 
 + c \int_1^t \tau ^{-3\ell/2}  \| \chi^{(n-1)} (\tau )\|_s \Big(\tau ^{3\ell/2}\| \chi^{(n-1)} (\tau )\|_s
+ \tau ^{3\ell/2}\|\Psi  (\tau )\|_s \Big)^2   \,d\tau \\
& \leq  & 
c   \varepsilon \|\chi _0 (x) \|_s 
 + c \int_1^t\tau ^{-3\ell} \Big( \tau ^{3\ell/2} \| \chi^{(n-1)} (\tau )\|_s \Big)\Big( \tau ^{3\ell/2}\| \chi^{(n-1)} (\tau )\|_s
+ \|\Psi_0\|_s \Big)^2   \,d\tau \,.
\end{eqnarray*}
Hence,
\begin{eqnarray} 
a_n(t) 
& \leq  & 
c   \varepsilon \|\chi _0  \|_s 
 + c \int_1^t\tau ^{-3\ell} a_{n-1}(\tau ) \big( 1+ A_{n-1}\big)^2   \,d\tau \nonumber \\
& \leq  & 
c   \varepsilon \|\chi _0  \|_s 
 + c \big( 1+ A_{n-1}\big)^2 \int_1^t \tau ^{-3\ell} a_{n-1}(\tau )   \,d\tau  \nonumber \\
 \label{4.30FLWR}
& \leq  & 
c   \varepsilon \|\chi _0  \|_s 
 + c \big( 1+ A_{n-1}\big)^2 \int_1^t \tau ^{-3\ell} a_{n}(\tau )   \,d\tau \,.
\end{eqnarray}
Denote
\[
y(t):= \int_1^t \tau ^{-3\ell} a_{n}(\tau )   \,d\tau \quad \mbox{\rm then} \quad  y'(t)=  t^{-3\ell}  a_{n}(t)  
\]
and according to (\ref{4.30FLWR}) we obtain
\begin{eqnarray} 
\label{4.111FLRW}
y'(t) 
& \leq  &
t^{-3\ell}c\varepsilon \| \chi _0  \|_{s}+ 2c\big( 1+ A_{n-1}\big)^2 t^{-3\ell}y(t)\,,  
\end{eqnarray}
which implies
\[
\frac{d}{d t}\left( \exp \left({\dsp - 2c\big( 1+ A_{n-1}\big)^2\int_1^t  \tau ^{-3\ell}\, d\tau } \right) y(t)\right)\\
\leq 
 \exp \left({\dsp - 2c\big( 1+ A_{n-1}\big)^2\int_1^t  \tau ^{-3\ell}\, d\tau } \right) t^{-3\ell}c\varepsilon \| \chi _0  \|_{s} .
\]
Since $y(1)=0$, after evaluating the integral  we obtain
\[
 \exp \left({\dsp - 2c\big( 1+ A_{n-1}\big)^2\int_1^t  \tau ^{-3\ell}\, d\tau } \right)  y(t)\leq 
c\varepsilon \| \chi _0  \|_{s}\int_1^t  \exp \left({\dsp - 2c\big( 1+ A_{n-1}\big)^2\int_1^{\tau_1}  \tau ^{-3\ell}\, d\tau } \right)\tau_1 ^{-3\ell}\,  d\tau_1 .
\]
and
\begin{eqnarray*} 
\hspace{-0.5cm} y(t)& \leq  &
c\varepsilon \| \chi _0  \|_{s} \exp \left({\dsp  2c\big( 1+ A_{n-1}\big)^2\int_1^t  \tau ^{-3\ell}\, d\tau } \right)   \int_1^t  \exp \left({\dsp - 2c\big( 1+ A_{n-1}\big)^2\int_1^{\tau_1}  \tau ^{-3\ell}\, d\tau } \right)\tau_1 ^{-3\ell}\,  d\tau_1 .
\end{eqnarray*}
Hence
\begin{eqnarray*} 
 y(t)
& \leq  &
c\varepsilon \| \chi _0  \|_{s} \int_1^t \exp \left({\dsp  2c\big( 1+ A_{n-1}\big)^2\int^t_{\tau_1}  \tau ^{-3\ell}\, d\tau } \right)\tau_1 ^{-3\ell}\,  d\tau_1 \,,
\end{eqnarray*}
and according to (\ref{4.111FLRW}) 
\begin{eqnarray*}
t^{-3\ell} a_{n}(t) 
& \leq  &
t^{-3\ell} c\varepsilon \| \chi _0  \|_{s}+ 2c\big( 1+ A_{n-1}\big)^2 t^{-3\ell} y(t)\\
& \leq  &
t^{-3\ell} c\varepsilon \| \chi _0  \|_{s}
+ 2c\big( 1+ A_{n-1}\big)^2 t^{-3\ell} c\varepsilon \| \chi _0  \|_{s} \int_1^t \exp \left({\dsp  2c\big( 1+ A_{n-1}\big)^2\int^t_{\tau_1}  \tau ^{-3\ell}\, d\tau } \right)\tau_1 ^{-3\ell}\,  d\tau_1,  
\end{eqnarray*}
that is, 
\[
a_{n}(t)  
\leq 
c\varepsilon \| \chi _0  \|_{s}\Bigg\{ 1 
+ 2c\big( 1+ A_{n-1}\big)^2\int_1^t \exp \left({\dsp  2c\big( 1+ A_{n-1}\big)^2\int^t_{\tau_1}  \tau ^{-3\ell}\, d\tau } \right)\tau_1 ^{-3\ell}\,  d\tau_1\Bigg\}\,,  
\]
and
\[ 
A_{n} 
\leq 
c\varepsilon \| \chi _0  \|_{s}\Bigg\{ 1
+ 2c\big( 1+ A_{n-1}\big)^2 \int_1^t \exp \left({\dsp  2c\big( 1+ A_{n-1}\big)^2\int^t_{\tau_1}  \tau ^{-3\ell}\, d\tau } \right)\tau_1 ^{-3\ell}\,  d\tau_1\Bigg\}\,.  
\]
On the other hand
\[
2c\big( 1+ A_{n-1}\big)^2 \int_1^t \exp \left({\dsp  2c\big( 1+ A_{n-1}\big)^2\int^t_{\tau_1}  \tau ^{-3\ell}\, d\tau } \right)\tau_1 ^{-3\ell}\,  d\tau_1 
 = 
 \exp \left({\dsp  2c\big( 1+ A_{n-1}\big)^2\int_{1}^t   \tau ^{-3\ell}\, d\tau }\right)  -1  
\]
and 
the condition (\ref{aintehgral})    lead  to
\[ 
A_{n} 
\leq 
c\varepsilon \| \chi _0  \|_{s} \exp \left\{ \dsp \big( 1+ A_{n-1}\big)^2 C\right\}\,,\quad C:=\frac{2c}{1-3\ell} \,.  
\]
Then, with some other constant $\widetilde C $ we obtain 
\begin{eqnarray} 
\label{4.32}
A_{n} 
& \leq  &
\widetilde C\varepsilon \| \chi _0  \|_{s}  \exp \left\{ \dsp  \widetilde CA_{n-1}\big( 1+ A_{n-1}\big)\right\}\,.
\end{eqnarray}
Let $ \varepsilon _0$ be such that  
\[
 {  2 \widetilde C^2\varepsilon_0   \| \chi _0  \|_{s} \big( 1+ 2\widetilde C\varepsilon_0   \| \chi _0  \|_{s} \big)  < \ln 2 }\,.
\]
If
\[
A_{n-1}\leq 2\widetilde C\varepsilon { \| \chi _0  \|_{s} }\quad \mbox{\rm and} \quad\varepsilon \leq \varepsilon _0 \,, 
\]
then due to (\ref{4.32})  we obtain
\begin{eqnarray*} 
A_{n} 
& \leq  &
{ \widetilde C\varepsilon \| \chi _0  \|_{s}  \exp \left\{ \dsp  \widetilde C2\widetilde C\varepsilon \| \chi _0  \|_{s} \big( 1+ 2\widetilde C\varepsilon \| \chi _0  \|_{s} \big)\right\} } \\
& \leq  &
2\widetilde C\varepsilon \| \chi _0  \|_{s}  \,.
\end{eqnarray*}
Thus, for given $s$ and for all $n \geq 1$ we have proved the estimate
\[
\sup_{t \in [1,\infty)}\, \sup_{{1\leq \tau \leq t\,,\,0\leq k\leq n}} \tau ^{3\ell /2}\| \chi^{(n)}  (\tau ) \|_{s}\leq 
2\widetilde C\varepsilon \| \chi _0  \|_{s} < \infty \quad \mbox{\rm  for all} \quad n \geq 1\,.
\]
The last estimate, (\ref{4.26FLRW}), and Sobolev inequality imply
\begin{eqnarray} 
\label{B47}
\sup_{ n=0,1,2\ldots } \,\,\sup_{x \in {\mathbb R}^3 ,\,\,t \in [1,\infty)}\left\{t^{3\ell /2} | \chi^{(n)}  (x,t) |\,,\, t^{3\ell /2} | \Psi (x,t) |   \right\} = r < \infty\,.
\end{eqnarray}
Hence, for all $(x,t) \in {\mathbb R}^3\times [1,\infty)$ we have
\begin{eqnarray}
\label{4.29} 
&  &
\Big|f_1(\chi^{(k-1)} (x,t); \Psi (x,t) ) -f_1(\chi^{(k-2)} (x,t); \Psi (x,t) )  \Big| \\
& \leq &
\Big|  \chi^{(k-1)} (x,t)  - \chi^{(k-2)} (x,t)   \Big| \sup_{{\xi , \eta  \in {\mathbb C}\,\,| \xi |, | \eta | \leq r}}   | \nabla_{\xi ,\eta}f_1 (\xi ,\eta ) |\,.\nonumber
\end{eqnarray}
Consider 
\begin{eqnarray*}  
&  &
\dsp\left(   i {\gamma }^0    \partial_0   +i t^{-\ell} \sum_{j=1,2,3}{\gamma }^j  \partial_{x_j} +i \frac{3\ell }{2 } t^{-1}   {\gamma }^0      -m{\mathbb I}_4 +i{\gamma }^0 V(x,t) \right)\Big( \chi^{(k)} -\chi^{(k-1)} \Big)\\
&  = &
f_1\left(\chi^{(k-1)} \,,\,\Psi  \right) - f_1\left(\chi^{(k-2)} \,,\,\Psi  \right)  \,,\quad k=1,2,\ldots\,.
\end{eqnarray*}
By Lemma~\ref{L1.2FLRW}
for the solution of the last equation taking into account (\ref{4.29}) and the initial values, 
one has 
\begin{eqnarray*}
&  &
t^{3\ell/2}\|  \chi^{(k)} (t)-\chi^{(k-1)}(t) \|_{L^2({\mathbb R}^3)} \\
& \leq  &
c\int_0^t  s^{3\ell/2}  \|f_1\left(\chi^{(k-1)}(x, \tau ) \,,\,\Psi (x, \tau )  \right) - f_1\left(\chi^{(k-2)}(x, \tau )  \,,\,\Psi (x, \tau )  \right) \|_{L^2({\mathbb R}^3)} \,ds\\ 
& \leq  &
c \Bigg( \sup_{{\xi , \eta  \in {\mathbb C},\,\| \xi |, | \eta | \leq r}}   |\nabla_{\xi ,\eta}f_1 (\xi ,\eta ) | \Bigg)\int_0^t s^{3\ell/2}    \|
 \chi^{(k-1)} (x,t)  - \chi^{(k-2)} (x,t)   \|_{L^2({\mathbb R}^3)} \,ds\,.
\end{eqnarray*}
It follows
\begin{eqnarray*}
&  &
t^{3\ell/2}\|  \chi^{(k)} (t)-\chi^{(k-1)}(t) \|_{L^2({\mathbb R}^3)} \\
& \leq  &
c^2 \Bigg( \sup_{{\xi , \eta  \in {\mathbb C},\,\,| \xi |, | \eta | \leq r}}   |\nabla_{\xi ,\eta}f_1 (\xi ,\eta ) | \Bigg)^2 \int_0^t \,ds \int_0^{s} s_1^{3\ell/2} \|
 \chi^{(k-2)} (x,s_1)  - \chi^{(k-3)} (x,s_1)   \|_{L^2({\mathbb R}^3)} \,ds_1 
\end{eqnarray*}
and, consequently,
\begin{eqnarray*} 
t^{3\ell/2}\|  \chi^{(k)} (t)-\chi^{(k-1)}(t) \|_{L^2({\mathbb R}^3)} 
\leq 
C\| \chi _0  \|_{L^2({\mathbb R}^3)}
\frac{(Ct)^k}{k!}, \quad \mbox{\rm for all} \quad k=1,2,\ldots\,. 
\end{eqnarray*}
Hence, the sequence $ t^{3\ell/2} \chi^{(k)} (t)$  converges  to some  $ t^{3\ell/2}\chi \in C^0([1,\infty); (L^2({\mathbb R}^3))^4)$, that is,
\[
\lim_{k \to \infty} t^{3\ell/2}  \chi^{(k)} (t) = t^{3\ell/2} \chi  (t) 
\]
uniformly on every compact subset of ${\mathbb R}^3 $. By (\ref{B47})
\[
\lim_{k \to \infty} f_1\left(\chi^{(k)} \,,\,\Psi  \right) = f_1\left(\chi \,,\,\Psi  \right)  \quad \mbox{\rm in} \quad C^0([1,\infty); (L^2({\mathbb R}^3))^4)\,.
\]
Thus, $\chi  $ solves (\ref{B34}) while $ \psi $ solves  (\ref{NDE_CP})  and
\begin{eqnarray}
\label{B48}
\sup_{t \in [1,\infty) }t^{3\ell/2} \| \psi   (t) \|_{s}  < \infty \quad \mbox{\rm for } \quad s \geq 3
\end{eqnarray}
implies 
\[
\sup_{t \in [1,\infty) }t^{3\ell/2}| \psi   (x,t)|_{s^\prime }  < \infty\quad \mbox{\rm for } \quad s^\prime  <s-3/2 \,.
\]
{ It follows that  } for any integer $s \geq 3$ we have
\begin{eqnarray}
\label{B49}
\left\| F\left(\ol{\psi  }^T\gamma ^0 \psi\,,\,\ol{\psi  }^T\gamma ^0\gamma ^5 \psi \right)\psi (t) \right\|_{(H_s({\mathbb R}^3))^4}
\leq C_s\left\| \psi (t)\right\|_{(H_s({\mathbb R}^3))^4}\,.
\end{eqnarray}
For every $s \geq 3 $ the local Cauchy problem for (\ref{NDE_CP}) is well posed  in $C^0([0,T_s);(H_s({\mathbb R}^3))^4)  $ for some $1<T_s$. According to (\ref{B49}) and
\begin{eqnarray*}
 \|\psi (x,t)\|_s   
& \leq   &
C_s t^{-3\ell/2}   + C_s  t^{-3\ell/2} \int_0^t \tau ^{3\ell/2} \|\psi (x,\tau )\|_s  \,  d\tau \,,\qquad t \in [1,T_s)\,,
\end{eqnarray*}
from the last inequality we conclude that $T_s=\infty$ and $\psi (t) \in (C_0^\infty({\mathbb R}^3))^4 $.  The equation (\ref{NDE_CP}) implies $\psi   \in (C^{1,\infty}_{t,x}([1,\infty) \times {\mathbb R}^3) )^4 $. Theorem is proved. \qed

\section{Asymptotics at infinity}
\label{S4}
\setcounter{equation}{0}

\begin{theorem}
\label{T4.1A}
Assume that all conditions of Theorem~\ref{T4.1} are fulfilled. Let $\psi =\psi (x,t) $ be a solution of the problem (\ref{NDE_CP})  
given by Theorem~\ref{T4.1}. 
Suppose also that  $F\left( \psi^* \gamma ^0 \psi\,,\,\psi^*\gamma ^0\gamma ^5 \psi \right)\psi$  is the Lipschitz continuous function with exponent $\alpha  >0 $ in the space $H_{(3)}({\mathbb R}^3) $, 
 where  $3\ell\alpha >2$. 

Then the limit
\begin{equation}
\label{LimInt}
\lim_{t \to  \infty}   \int_1^{ t} S(1,\tau ) F\left(\psi^* (x,\tau )\gamma ^0 \psi(x,\tau )\,,\,\psi^* (x,\tau )\gamma ^0\gamma ^5 \psi (x,\tau )\right)\psi (x,\tau )\,d \tau   
\end{equation}
exists in the space  $(H_{(3)}({\mathbb R}^3))^4 $.
Furthermore, the solution $\widetilde \psi (x,t)$ of the Cauchy problem for the free Dirac equation 
\begin{eqnarray*} 
\cases{ \dsp\left( {\mathscr D}_{FLRW}(t,\partial_t,\partial_x)  + {\gamma }^0 V(x,t) \right)\widetilde \psi  
{ (x,t)}   =0  \,,\cr
\widetilde \psi (x,1)=\widetilde \psi (x)\,,
}
\end{eqnarray*}
where 
\[ 
\widetilde \psi (x) 
  =   
\Psi _0 (x)+\varepsilon \chi _0  (x)
+\lim_{t \to  \infty}   \int_1^{ t} S(1,\tau ) F\left(\psi^* (x,\tau )\gamma ^0 \psi(x,\tau )\,,\,\psi^* (x,\tau )\gamma ^0\gamma ^5 \psi (x,\tau )\right)\psi (x,\tau )\,d \tau\,,
\] 
satisfies  
\[
\lim_{t \to + \infty}    \left\| \psi (x,t)-  \widetilde \psi (x,t) \right\|_{(H_{(3)}({\mathbb R}^3))^4}=0\,.
\]

Moreover, if $V(x,t)=0 $, then
\[
\widetilde \psi  (x,t) 
  =   {\mathscr{D}}^{co} (t,\partial _t,\partial _x)\gamma^0\left (
   \begin{array}{ccccc}
  {\cal K}_1(x,t,D_x;m;1 ) {\mathbb I}_2& {\mathbb O}_2   \\
   {\mathbb O}_2&  {\cal K}_1(x,t,D_x;-m;1 ){\mathbb I}_2   \\ 
   \end{array}
   \right) [\widetilde \psi (x)]  (x,t),\quad t>1\,.
\]
\end{theorem}
\medskip

\noindent
{\bf Proof.} It is enough to prove the convergence in the space $(H_{(3)}({\mathbb R}^3))^4$. 
The solution $ \psi $ can be written as follows (\ref{B33})
\begin{eqnarray*} 
\psi =\Psi +\chi \,,
\end{eqnarray*}
where the function $\Psi=\Psi (x,t) $ solves   (\ref{B31}) 
\[  
\cases{ \dsp\left(  {\mathscr D}_{FLRW}(t,\partial_t,\partial_x)  + {\gamma }^0 V(x,t) \right)\Psi  =0  \,,\cr
\Psi (x,1)=\Psi _0 (x)\,,
}
\]
while the function $\chi =\chi  (x,t) $ solves   (\ref{B34})
\[
\cases{ \dsp\left(   {\mathscr D}_{FLRW}(t,\partial_t,\partial_x)  + {\gamma }^0 V(x,t) \right)\chi   =f_1\left(\chi \,,\,\Psi  \right)   \,,\cr
\chi  (x,1)=\varepsilon \chi _0 (x)\,,
}
\]
where
\[
|f_1\left(\chi \,,\,\Psi  \right)| =O\Big( |\chi |\big(|\chi |+|\Psi  |\big)^2\Big)\,.
\]
Since $F\left( \psi^* \gamma ^0 \psi\,,\,\psi^*\gamma ^0\gamma ^5 \psi \right)\psi$  is a Lipschitz continuous function with exponent $\alpha  >0 $ and according to 
(\ref{B48}) and (\ref{BackFLRW}) of Lemma~\ref{L1.2FLRW}, 
\begin{eqnarray*}
\left\| F\left( \psi^*\gamma ^0 \psi\,,\,\psi^*\gamma ^0\gamma ^5\psi \right)\psi (\tau )\right\|_{(H_{(3)}({\mathbb R}^3))^4}
 & \leq  &
C_s\left\| \psi  (\tau ) \right\|_{(H_{(3)}({\mathbb R}^3))^4}^{1+\alpha  }\\
 & \leq  &
C_s\tau^{-\frac{3}{2}\ell(1+\alpha )}\,.
\end{eqnarray*}
At the same time
\begin{eqnarray*}
\left\| S(1,\tau ) F\left(\psi^*\gamma ^0 \psi\,,\,\psi^*\gamma ^0\gamma ^5\psi \right)\psi (\tau )\right\|_{(H_{(3)}({\mathbb R}^3))^4}
& \leq  &
C_s\tau^{\frac{3}{2}\ell}\left\|  F\left(\psi^*\gamma ^0 \psi\,,\,\psi^*\gamma ^0\gamma ^5\psi \right)\psi (\tau )\right\|_{(H_{(3)}({\mathbb R}^3))^4}\\
& \leq  &
C_s\tau ^{-\frac{3}{2}\ell\alpha }\,.
\end{eqnarray*}
The last estimate due to the condition $3\ell\alpha >2$ implies a convergence of  (\ref{LimInt}) in $(H_{(3)}({\mathbb R}^3))^4 $.
Theorem is proved. 
\qed

{ We skip a statement on the continuity of the scattering operator similar to one made in Theorem~\ref{T1.2}, since it can be easily recovered from the calculations written above, while  a proof of such statement is really  cumbersome and does not contain any new ideas.}

\medskip

\section{Nonexistence of global in time solution}
\label{S5}
\setcounter{equation}{0}

Consider the semilinear  Dirac equation 
\[
  \dsp 
\left(  i {\gamma }^0    \partial_0   +i t^{-\ell}\sum_{j=1,2,3}{\gamma }^j  \partial_{j} +i\frac{3\ell }{2 } t^{-1}    {\gamma }^0     -m t^{-1}{\mathbb I}_4 +{\gamma }^0V(x,t)\right)\psi=iG(\psi ){\gamma }^0\psi  \,,
\]
where  $m \in \C$,  and the matrix-valued term  $ G(\psi )$ commutes with ${\gamma }^0 $,
\[
{\gamma }^0 G(\psi )=G(\psi ){\gamma }^0 . 
\]
The equation can  also be written in the equivalent form of the following symmetric hyperbolic system 
\[
  \dsp 
\left(       \partial_0   +   t^{-\ell}\sum_{j=1,2,3}{\alpha }^j  \partial_j+  \frac{3\ell }{2 } t^{-1}{\mathbb I}_4     +i mt^{-1} {\gamma }^0 -iV (x,t)\right)\psi=G(\psi ) \psi  \,.
\]
Consider the energy integral
\[
E(t)
  =  
\int_{{\mathbb R}^3}|\psi (x,t)|^2\,dx  .
\] 
\begin{lemma}
\label{L1}
Let   $m \in \C$  and $V^*(x,t)=V(x,t)$. Then
\[
\frac{d}{d t} E(t)
=
  \int_{{\mathbb R}^3}\left(  2\Re ( { \ol{G_{jk}(\psi )   \psi_k (x,t) }\psi_j (x,t) })  - 3\ell t^{-1}|\psi (x,t)|^2  + 2t^{-1}(\Im (m)) \psi^* (x,t)  {\gamma }^0 \psi (x,t)    \right) dx.
\]
\end{lemma} 
\medskip

\ndt
{\bf Proof.} The proof is similar to the proof of Lemma~\ref{L1.2FLRW} and we skip it. 
\qed

The next theorem gives the blowup result for the Dirac equation in the  universe with  { the scale function $a(t)=t^\ell$, $\ell \in  \mathbb R $. }

\begin{theorem}
\label{TBUH+}
Consider the Cauchy problem 
\begin{equation}
\label{CPHP}
\cases{
  \dsp 
\left(  i {\gamma }^0    \partial_0   +i t^{-\ell}\sum_{j=1,2,3}{\gamma }^j  \partial_{j} +i \frac{3\ell }{2 } t^{-1}    {\gamma }^0     -m t^{-1}{\mathbb I}_4  +{\gamma }^0V(x,t)\right)\psi=G(\psi )i{\gamma }^0\psi  \,,\cr
\psi (x,1) =\psi_0 (x )
}
\end{equation}
with $ \psi_0 (x )$ such that $
\mbox{\rm supp\,} \psi_0 (x ) \subseteq  \{x \in {\mathbb R}^3\,|\,|x| \leq R\}$.
Assume that
\begin{eqnarray*}
&  &
 G(\zeta ) =O(|\zeta |),  \quad  {\gamma }^0 G(\zeta )=G(\zeta ){\gamma }^0 , \hspace*{0.5cm} {  \Re (\zeta^*G(\zeta )^*\zeta  )}\geq c_0|\zeta |^{2+\alpha },\quad { c_0>0},\quad \alpha >0\,.
\end{eqnarray*}
Then for  arbitrary size initial data there is no global solution  $\psi $ of (\ref{CPHP}) that obeys the finite propagation speed property
\begin{eqnarray*}
&  &
\mbox{\rm supp\,} \psi (x,t) \subseteq 
\left\{x \in {\mathbb R}^3\,\Big|  \,|x| \leq R +{\mathcal A}(t)  \right\}, \quad {\mathcal A}(t):=\int_1^t \frac{1}{a(\tau )}\,d\tau\,, 
\end{eqnarray*}
 in the following cases:  
 \begin{itemize}
\item  If $\ell<1 $ and $1\geq \frac{3\alpha}{2}+\alpha|\Im (m)| $, 
\item   If  $1<\ell  $ and $1\geq \frac{3\alpha}{2}\ell+\alpha|\Im (m)| $, 
\item  If $ \ell=1 $ and $1\geq  \frac{3\alpha}{2}+\alpha|\Im (m)| $.   
\end{itemize}

Furthermore, in the following cases for the large  initial data there is no global solution of (\ref{CPHP}) that obeys the finite propagation speed property:
 \begin{itemize}
\item  If  $\ell<1 $ and $1<\frac{3\alpha}{2}+\alpha|\Im (m)| $. 
\item  If $1<\ell  $ and $1<\frac{3\alpha}{2}\ell+\alpha|\Im (m)| $. 
\item  If $ \ell=1 $ and $1< \frac{3\alpha}{2}+\alpha|\Im (m)| $. 
\end{itemize}
 The solution blows up no later than  time $T_{bu}$ depending on  $c_0$, $\alpha $, $m$, $\psi _0$, and $\ell$  such that 
\begin{eqnarray} 
\label{Tbu} 
& & 
{ \left( E(1)  \right)^{ -\frac{\alpha}{2} } =  \frac{\alpha}{2}  c_0\int_1^{T_{bu} }   \left( R +{\mathcal A}(t) \right)^{-\frac{3\alpha}{2} }    t^{-\ell \frac{3\alpha}{2}-\alpha |\Im (m)|}\,dt }
\end{eqnarray}
and  
\[
 \lim_{t \nearrow T_{bu}}\int_{{\mathbb R}^3}|\psi (x,t )|^{2} \,dx =\infty\,,
 \]
provided that $E(1)$ is sufficiently large.
\end{theorem}
\medskip

\ndt
{\bf Proof.} According to 
Lemma~\ref{L1} 
if    $m \in \C$, then
\begin{eqnarray}
\frac{d}{d t} E(t)
& \geq  &
  \int_{{\mathbb R}^3}\left( c_0 |\psi (x,t)|^{2+\alpha } - 3\ell  t^{-1} |\psi (x,t)|^2  + 2t^{-1}(\Im (m)) \psi^* (x,t)  {\gamma }^0 \psi (x,t)     \right) \,dx \nonumber \\
& \geq  &
  c_0\int_{{\mathbb R}^3}|\psi (x,t)|^{2+\alpha } \,dx- \int_{{\mathbb R}^3} \left( 3\ell  t^{-1}  |\psi (x,t)|^2 +  2t^{-1}|\Im (m)|  |\psi (x,t)|^2   \right) \,dx \nonumber \\
  \label{5.3b}
& \geq  &
c_0  \int_{{\mathbb R}^3}|\psi (x,t)|^{2+\alpha } \,dx-\left( 3\ell  t^{-1} + 2t^{-1}|\Im (m)| \right)\int_{{\mathbb R}^3} |\psi (x,t)|^2  \,dx\,.
\end{eqnarray}
Further, if the solution obeys the finite speed propagation property, we obtain
\begin{eqnarray*}
\int_{{\mathbb R}^3}  |\psi (x,t)|^2  \,dx
& \leq  &
\left( \int_{\mbox{\rm supp\,} \psi}   |\psi (x,t)|^{2+\alpha }  \,dx \right)^{2/(2+\alpha )}\left( R +{\mathcal A}(t)  \right)^{3\alpha /(2+\alpha )}  \,.
\end{eqnarray*}
It follows
\[
\left( R +{\mathcal A}(t)  \right)^{-3\alpha /2} \left( \int_{{\mathbb R}^3}  |\psi (x,t)|^2  \,dx\right)^{(2+\alpha )/2}
 \leq  
\int_{\mbox{\rm supp\,} \psi}   |\psi (x,t)|^{2+\alpha }  \,dx \,.
\]
Then by (\ref{5.3b})
\begin{eqnarray}
\frac{d}{d t} E(t)
\label{5.3}
& \geq  &
c_0  \left( R +{\mathcal A}(t) \right)^{-3\alpha /2}   E(t) ^{(2+\alpha )/2} 
-\left( 3\ell  t^{-1} + 2t^{-1}|\Im (m)| \right)E(t)\,.  
\end{eqnarray}
Denote
\begin{eqnarray*}
K(t)
& :=  &
c_0  \left( R +{\mathcal A}(t)  \right)^{-3\alpha /2}\,, \\
A(t) & :=  & \int_1^t \left(  3\ell \tau ^{-1} + 2\tau ^{-1}|\Im (m)|\right)\,d \tau=
\ln a(t)^3 - \ln a(1)^3 + 2|\Im (m)|\ln (t)   ,
\end{eqnarray*}
then the inequality  (\ref{5.3}) reads 
\begin{eqnarray*}
\frac{d}{d t} E(t)
& \geq  & 
K(t) E(t)^{(2+\alpha )/2} 
-E(t)\frac{d}{d t}A(t)\,.
\end{eqnarray*}
Hence
\begin{eqnarray*}
\frac{d}{d t} \left( E(t)e^{A(t)} \right) 
& \geq  & 
K(t)   e^{-A(t)\alpha /2}     \left( E(t)e^{A(t)} \right)^{(2+\alpha )/2}  \,.
\end{eqnarray*}
For the function
\begin{eqnarray*}
 F(t)
& :=  & 
  E(t)e^{A(t)}    
\end{eqnarray*}
the inequality leads to
\begin{eqnarray*}
\frac{d}{d t} F(t) 
& \geq  & 
K(t)   e^{-A(t)\alpha    /2}     F(t) ^{(2+\alpha )/2}  
\end{eqnarray*}
or
\begin{eqnarray*}
 \frac{d}{d t} F(t) ^{ -\alpha /2} 
& \leq   & 
-\frac{\alpha}{2} K(t)   e^{-A(t)\frac{\alpha}{2}   }  \,. 
\end{eqnarray*}
After integration we obtain
\begin{eqnarray*}
   F(t) ^{ -\frac{\alpha}{2}}   
& \leq   & 
F(1) ^{ -\frac{\alpha}{2}} -\frac{\alpha}{2}  c_0\int_1^t   \left( R +{\mathcal A}(s) \right)^{-3\frac{\alpha}{2}}   e^{-A(s)\frac{\alpha}{2} } \,ds  \,.   
\end{eqnarray*} 
Furthermore 
\begin{eqnarray*} 
\int_1^t   \left( R +{\mathcal A}(s)\right)^{-3\frac{\alpha}{2}}   e^{-A(s)\frac{\alpha}{2} } \,ds  
& = & 
\int_1^t   \left( R +{\mathcal A}(s) \right)^{-3\frac{\alpha}{2}} \left( \frac{a(1) }{a(s) }\right)^{ \frac{3\alpha}{2}}  s^{-\alpha|\Im (m)|}  \,ds \\
& = & 
J(t,m,R, \alpha ,\ell)\,,  
\end{eqnarray*}
where 
\begin{eqnarray*}
J(t,m,R, \alpha ,\ell)
& := & 
 \int_1^t   \left( R +{\mathcal A}(s) \right)^{-3\alpha /2}  s^{- \ell\frac{3\alpha}{2}-\alpha|\Im (m)|}  \,ds\,.   
\end{eqnarray*}
\noindent
In the   case of $\ell<1 $ and $1= \frac{3\alpha}{2}+\alpha|\Im (m)| $ we obtain
 \begin{eqnarray*}
J(t,m,R, \alpha ,\ell) & \approx  &
C + \int_1^t   \left(   s^{1-\ell} \right)^{-\frac{3\alpha}{2}}  s^{- \ell\frac{3\alpha}{2}-\alpha|\Im (m)|}  \,ds= C + \int_1^t   s^{-1}  \,ds = \ln(t) \to \infty \quad \mbox{for large}\,\, t.   
\end{eqnarray*}

\noindent  
In the   case of $\ell<1 $ and $1> \frac{3\alpha}{2}+\alpha|\Im (m)| $ we obtain
 \begin{eqnarray*}
J(t,m,R, \alpha ,\ell)& \approx  &
C + \int_1^t   \left(   s^{1-\ell} \right)^{-\frac{3\alpha}{2}}  s^{- \ell\frac{3\alpha}{2}-\alpha|\Im (m)|}  \,ds\approx  C +   t^{1-\frac{3\alpha}{2}-\alpha|\Im (m)|} \to \infty \quad \mbox{for large}\quad t\,.   
\end{eqnarray*}  
Thus, in these two cases the solution blows up for { initial data with  arbitrary size}.\\

\noindent
In the   case  of $\ell<1 $ and $1<\frac{3\alpha}{2}+\alpha|\Im (m)| $ we obtain
 \begin{eqnarray*}
J(t,m,R, \alpha ,\ell)& \approx  &
C + \int_1^t   \left(   s^{1-\ell} \right)^{-\frac{3\alpha}{2}}  s^{- \ell\frac{3\alpha}{2}-\alpha|\Im (m)|}  \,ds\leq 2C   \quad \mbox{for large}\quad t\,.   
\end{eqnarray*}  
Thus, the solution blows up for the large  initial data. The solution blows up no later than  time $T_{bu}$ such that (\ref{Tbu})  holds and  $E(1)$ is sufficiently large. \\

\noindent
In the   case of $1<\ell  $ and $1=\frac{3\alpha}{2}\ell+\alpha|\Im (m)| $ we obtain
\begin{eqnarray*}
J(t,m,R, \alpha ,\ell)& \approx  &  
C + \int_1^t   \left( R +\frac{1}{\ell-1} s^{1-\ell}\right)^{-3\alpha /2}  s^{- \ell\frac{3\alpha}{2}-\alpha|\Im (m)|}    \,ds\\
& \approx  &  
C +  2R ^{-3\alpha /2}  \ln (t) \to \infty \quad \mbox{for large}\quad t\,.   
\end{eqnarray*} 

\noindent
In the   case of $1<\ell  $ and $1>\frac{3\alpha}{2}\ell+\alpha|\Im (m)| $ we obtain
\begin{eqnarray*}
J(t,m,R, \alpha ,\ell)& \approx  &  
C + \int_1^t   \left( R +\frac{1}{\ell-1}s^{1-\ell} \right)^{-3\alpha /2}  s^{- \ell\frac{3\alpha}{2}-\alpha|\Im (m)|}    \,ds\\
& \approx  &  
C +  2R ^{-3\alpha /2} \frac{1}{1- \ell\frac{3\alpha}{2}-\alpha|\Im (m)|}  t^{1- \ell\frac{3\alpha}{2}-\alpha|\Im (m)|}  \to \infty \quad \mbox{for large}\quad t\,.   
\end{eqnarray*}
Thus, in the last two cases the solution blows up for { initial data with  arbitrary size}.\\

\noindent
In the   case of $1<\ell  $ and $1<\frac{3\alpha}{2}\ell+\alpha|\Im (m)| $ we obtain
\begin{eqnarray*}
J(t,m,R, \alpha ,\ell)
& \approx  &  
C +   2R^{-3\alpha /2} \frac{1}{1- \ell\frac{3\alpha}{2}-\alpha|\Im (m)|}  t^{1- \ell\frac{3\alpha}{2}-\alpha|\Im (m)|}  \leq 2C  \quad \mbox{for large}\quad t\,.   
\end{eqnarray*}
Thus, the solution blows up for the large  initial data. The solution blows up no later than  time $T_{bu}$ such that (\ref{Tbu}) holds and 
 $E(1)$ is sufficiently large. \\

\noindent
In the   case of $ \ell=1 $ and $1= \frac{3\alpha}{2}+\alpha|\Im (m)| $ we obtain (for $\alpha \not=  2/3$)
 \begin{eqnarray*}
J(t,m,R, \alpha ,\ell)& \approx  &
C + \int_2^t   \left(   \ln (s) \right)^{-\frac{3\alpha}{2}}   s^{- \ell\frac{3\alpha}{2}-\alpha|\Im (m)|}  \,ds= C +\int_2^t  \left(   \ln (s) \right)^{-\frac{3\alpha}{2}}  s^{-1}  \,ds\\
&  = &
C +\frac{1}{1-\frac{3\alpha}{2}} \left(   (\ln (t)) ^{1-\frac{3\alpha}{2}}-(\ln (2))^{1-\frac{3\alpha}{2}}\right) \quad \mbox{for large}\quad t\,.   
\end{eqnarray*}  
Hence, there is no global solution for { initial data with  arbitrary size}. \\

\noindent
In the   case of $ \ell=1 $ and $1< \frac{3\alpha}{2}+\alpha|\Im (m)| $ we obtain
\begin{eqnarray*}
J(t,m,R, \alpha ,\ell)
& \approx  &
C + \int_1^t   \left(   \ln (s) \right)^{-\frac{3\alpha}{2}} s^{- \frac{3\alpha}{2}-\alpha|\Im (m)|}  \,ds\leq 2C   \quad \mbox{for large}\quad t\,.   
\end{eqnarray*}  
There is no global solution for large data. The solution blows up no later than  time $T_{bu}$ such that (\ref{Tbu}) 
 holds and  $E(1)$ is sufficiently large.\\

\noindent
In the   case of $ \ell=1 $ and $1> \frac{3\alpha}{2}+\alpha|\Im (m)| $ we obtain
\begin{eqnarray*}
J(t,m,R, \alpha ,\ell)
& \approx  &
C + \int_1^t   \left(   \ln (s) \right)^{-\frac{3\alpha}{2}} s^{- \frac{3\alpha}{2}-\alpha|\Im (m)|}  \,ds \to \infty  \quad \mbox{for large}\quad t\,.   
\end{eqnarray*}  
There is no global solution  for { initial data with  arbitrary size}. \\
Theorem~\ref{TBUH+} is proved. \qed

\begin{corollary} 
For the  matter dominated universe  with  $\ell =2/3$ and for the  radiation dominated universe  with  $\ell =1/2$  
 and $1\geq \frac{3\alpha}{2}+\alpha|\Im (m)| $,  for  {\bf arbitrary size} initial data there is no global solution  $\psi $ of (\ref{CPHP}) that obeys the finite propagation speed property. If   $1<\frac{3\alpha}{2}+\alpha|\Im (m)| $, then  for the {\bf large  initial data} there is no global solution of (\ref{CPHP}) that obeys the finite propagation speed property.\\
\end{corollary}

\end{document}